\input amstex
\loadbold
\loadeurm
\loadeusm
\openup.8\jot \magnification=1200

\def\even{{\text{even}}}
\def\M{{\Cal M}}
\def\J{{\Cal J}}
\def\G{{\Cal G}}

\def\cls{{\text{ cls}}}
\def\Z{\Bbb Z}
\def\A{\Bbb A}
\def\Q{\Bbb Q}
\def\disc{{\text{ disc}}}

\def\rank{{\text {rank}}}
\def\mod{{\text{mod\ }}}

\def\sgn{{\text{sgn}}}

\def\H{{\Bbb H}}
\def\rad{{\text {rad}}}
\def\e{{\text{e}}}

\def\h{{\it h}}
\def\R{{\Bbb R}}
\def\C{{\Bbb C}}

\def\mult{{\text{mult}}}
\def\diag{\text{diag}}

\def\gen{{\text{gen}}}

\def\cls{{\text{cls}}}
\def\F{{\Bbb F}}

\def\h{{\Cal H}}
\def\sym{{\text{sym}}}

\documentstyle{amsppt}
\pageheight{7.7in}
\vcorrection{-0.05in}
\topmatter
\pageno 1
\title 
A formula for the action of Hecke operators on half-integral weight Siegel
modular forms and applications
\endtitle

\rightheadtext{Hecke operators on half-integral weight Siegel forms}

\author Lynne H. Walling\endauthor
\subjclassyear{2010}
\subjclass 11F46, 11F60, 11F37, 11F27, 11F30\endsubjclass
\keywords  Siegel modular forms, 
Hecke operators, half-integral weight, theta series, Fourier coefficients\endkeywords
\address L.H. Walling, Department of Mathematics, University Walk,
University of Bristol, Bristol BS8 1TW England (telephone) (0) 117 331 5245 (fax) (0) 117
928 7999 \endaddress
\email l.walling\@bristol.ac.uk\endemail

\abstract 
We introduce
an alternate set of generators for the Hecka algebra,
and give an explicit formula for the action of these operators
on Fourier coefficients.  With this,
we compute the eigenvalues of Hecke operators acting on average
Siegel theta series with half-integral 
weight (provided the prime associated to the operators does
not divide the level of the theta series). 
Next, we bound the eigenvalues of these operators in terms
of bounds on Fourier coefficients.
Then we show that the
half-integral weight Kitaoka
subspace is stable under all Hecke operators.  Finally, we observe
that an obvious
isomorphism between Siegel modular forms of weight $k+1/2$ and ``even''
Jacobi modular forms of weight $k+1$ is Hecke-invariant (here the level
and character are arbitrary).
\endabstract

\endtopmatter

\document

\head{\bf \S0. Introduction 
}\endhead

Modular forms are of great interest in modern number theory, 
one reason being that
their Fourier coefficients carry number theoretic information.
A theta series attached to a positive definite quadratic form $Q$
is an example of this:  Given
$L=\Z v_1\oplus\cdots\oplus\Z v_m$ scaled
 so that $Q(L)\subseteq2\Z$, set
$$\theta(L;\tau)=\sum_{\ell\in L}e^{\pi i Q(\ell)\tau},$$
where $\tau$ lies in the complex upper half-plane;
it is well-known that $\theta(L;\tau)$ is a modular form of weight
$m/2$.
As a Fourier
series, $\theta(L;\tau)=\sum_{t\in\Z}r(L,2t)\,e^{2\pi it\tau}$ where
the Fourier coefficients are the representation numbers
$$r(L,2t)=\#\{\ell\in L: \ Q(\ell)=2t\ \}.$$
For each prime $p$ we have a Hecke operator $T(p)$ acting on the space
of modular forms;
Hecke operators help us study Fourier coefficients of modular forms,
as the space of modular forms has a basis of Hecke-eigenforms, and the
Fourier coefficients of a Hecke-eigenform satisfy arithmetic
relations.

Siegel was interested in generalised representation numbers that tell us the
number of rank $n$ sublattices $\Lambda$ of $L$ on which $Q$ restricts to $T$, $T$ any other
quadratic form.  For this he introduced generalised theta series, giving us the
first examples of Siegel modular forms.  
These modular forms have Fourier series
supported on symmetric, $n\times n$, even integral matrices (so the diagonal
entries of these matrices are even).
We again have Hecke operators, but for each prime $p$, 
we now have operators
$T(p)$ and $T_j(p^2)$ ($1\le j\le n$) generating the local Hecke algebra.
For integral weight, the action of $T(p)$ on Fourier coefficients is given in
[9], and the action of $T_j(p^2)$ on Fourier coefficients is given in [6].
In this paper we turn our attention to the half-integral weight case.

Hecke operators on half-integral weight Siegel modular forms, particularly
Siegel theta series, have been
studied before, noteably by Zhuravl\"ev in [17].  There the author
studies the image of the Hecke algebra under the Rallis map,
and subsequently the action on theta series with spherical harmonics.  He
develops formulas describing generalised Brandt matrices giving the action of
the Hecke operators, and gives an Euler expansion, in terms of these Brandt
matrices, for symmetric Dirichlet 
series built from the Siegel theta series.  He also obtains conditions for
linear dependence of theta series.  

Here we begin by considering
the standard generators of the Hecke algebra;
following [6], we compute their action on Fourier coefficients (see
Proposition 2.1 and Theorem 2.4).
However, the formulas for this action involve generalised ``twisted''
Gauss sums; while the values of these Gauss sums is explicitly known (proved in [11],
repeated in Proposition 1.4),
they are sufficiently complicated so that computations with these operators
are cumbersome.

Here we introduce an alternate set of generators $\widetilde T_j(p^2)$ for
the Hecke algebra; taking a very direct approach, we obtain an explicit
formula for the action of these operators on Fourier coefficients
(Theorem 3.3).  This formula easily leads to
several applications:  First, we extend [13], [14] to reprove Siegel's result
that, with $L$ an odd rank lattice equipped with a positive definite quadratic
form, the average theta series $\theta^{(n)}(\gen L)$ is an eigenform for the
Hecke operators attached to primes not dividing the level.  Further, we explicitly
compute the eigenvalues of the average theta series (Theorem 4.3), and we
identify Hecke operators that annihilate
the (unaveraged) theta series (Theorem 4.5).  Next, we bound the eigenvalues of these
Hecke operators in terms of bounds on Fourier coefficients (Theorem 5.1).
Then we give a quick proof that the ``Kitaoka space'' of
half-integral weight is invariant under all Hecke operators,
where the Kitaoka space
consists of those forms whose Fourier coefficients depend only on the genus
of the parameter (Theorem 6.2).  Finally, we observe that the formula for the Fourier
coefficients of $f|\widetilde T_j(p^2)$ is virtually identical to the formula in [15]
for the Fourier coefficients of $F|\widetilde T_j(p^2)$ where $F$ is a Jacobi modular
form of index 1; from this we see that with $\theta^{(n)}(\tau,Z)$ the classical
Jacobi theta series, $f(\tau)\mapsto f(\tau)\theta^{(n)}(\tau,Z)$ is a Hecke-invariant
isomorphism from the space of modular forms of level $N$,
character $\chi$, weight $k+1/2$ onto the
space of ``even'' Jacobi modular forms of level $N$, character $\chi'$, weight $k+1$, and index 1
(Theorem 7.4).
(We say a Jacobi modular form is even if it is supported only on pairs $(T,R)$ where
$R$, which multiplies $Z$ in the exponential, is even.  The relation between the characters
$\chi$ and $\chi'$ is stated and explained in Theorem 7.4.)

Note that Ibukiyama has a result related to the last application of
our formula.  In [7], without explicitly 
knowing the action of Hecke operators on
Fourier coefficients, Ibukiyama shows there is a Hecke-invariant
isomorphism between a Kohnen-type subspace
$\M^+_{k+1/2}(\Gamma^{(n)}_0(4),1)$ and $J_{k+1,1}(\Gamma_0^{(n,1)}(1),1)$.
Presuming there is a way to define a Kohnen-type subspace
$\M_{k+1/2}^+(\Gamma_0^{(n)}(4N),\chi)$ for arbitrary level $N\in\Z_+$
and character $\chi$, and
that there is a Hecke-invariant isomorphism
$\nu:\M_{k+1/2}^+(\Gamma_0^{(n)}(4N,\chi)\to
J_{k+1,1}(\Gamma_0^{(n,1)}(N),\chi'),$ we should have the following:  Let
$\eta$ denote our isomorphism from
$J^{\text{even}}_{k+1,1}(\Gamma_0^{(n,1)}(4N),\chi')$ onto
$\M_{k+1/2}(\Gamma_0^{(n)}(4N),\chi)$, and let $B_4$ be the ``shift
operator'' mapping $f(\tau)$ to $f(4\tau)$.  Then $B_4$ should map
$\M_{k+1/2}(\Gamma_0^{(n)}(4N),\chi)$ into 
$\M^+_{k+1/2}(\Gamma_0^{(n)}(16N),\chi)$, and $\nu\circ
B_4\circ\eta$ should be the identity map.

The author thanks Cris Poor, David Yuen, Jurg Kramer,
Martin Raum, and Ricard Hill for helpful and interesting
conversations. 

\bigskip
\head{\bf \S1. Preliminaries}
\endhead
\smallskip

Here we review the basic definitions, terminology, and results we rely on in 
the rest of the paper.  To read more about the basic theory of Siegel modular
forms, see for instance, [1], [2], [4]; to read more about the
basic theory of quadratic forms, see for instance, [5], [10].

Given odd $m\in\Z_+$ and a
lattice $L=\Z v_1\oplus\cdots\oplus\Z v_m$
 equipped with a positive definite quadratic form $Q$, and given
$n\in\Z$ with $n>1$, Siegel's generalised theta series is
$$\theta^{(n)}(L;\tau)=\sum_{G\in\Z^{m,n}}\e\{\,^tGQG\tau\}$$
where $\tau\in\h_{(n)}=\{X+iY:\ X,Y\in\R^{n,n}_{\sym},\ Y>0\ \}$
and $\e\{*\}=\exp(\pi iTr(*));$ here we have associated $Q$ with a symmetric
$m\times m$ matrix relative to the given basis for $L$.
We assume that $Q$ has been scaled to be even integral, meaning that
$^t\underline xQ\underline x\in 2\Z$ for $\underline x\in\Z^{m,1}$ (and hence
as a matrix, $Q$ is integral with even diagonal entries).  
The level of $L$ is the smallest positive integer $N$ so that
$NQ^{-1}$ is even integral.  (Note that since $L$ has odd rank, $4|N$;
see Theorem 8.9 of [5].)
As stated precisely in Theorem 1.2, and proved, for instance, in \S1
of [1], $\theta^{(n)}(L;\tau)$ transforms under the congruence
subgroup
$$\Gamma_0^{(n)}(N)=\left\{\pmatrix A&B\\C&D\endpmatrix\in Sp_n(\Z):\
C\equiv0\ (\mod N)\ \right\}$$
(here $Sp_n(\Z)$ is the symplectic group over $\Z$; in our notation,
the elements of $Sp_n(\Z)$ are $2n\times 2n$ matrices).  Note that
elements $\gamma=\pmatrix A&B\\C&D\endpmatrix\in GSp_n^+(\Q)$ act on
Siegel's upper half-space $\h_{(n)}$
by $\gamma \tau=(A\tau+B)(C\tau+D)^{-1}.$

We will need the following result, which is proved, for instance,
in Lemma 1.3.15 of [1].

\proclaim{Theorem 1.1 (Inversion Formula)} Let $L$ be a rank $m$ lattice, $m$ odd.
Define the dual of $L$ to be
$$L^{\#}=\{w\in\Q L:\ B_Q(w,L)\subseteq\Z\ \},$$
where $B_Q$ denotes the symmetric bilinear form associated to $Q$ via the
relation $Q(x+y)=Q(x)+Q(y)+2B_Q(x,y).$  Also, for $G_0\in\Q^{m,n}$, define
the inhomogeneous theta series
$$\theta^{(n)}(L,G_0;\tau)=\sum_{G\in\Z^{m,n}}\e\{Q[G+G_0]\tau\},$$
where $Q[E]=^tEQE$.  Then
$$\align
&\theta^{(n)}(L,G_0;\tau)\\
&\quad
=(\det Q)^{-n/2}(\det(-\tau))^{-m/2}\,\sum_{G\in\Z^{m,n}}\e\{-Q^{-1}[G]\tau^{-1}
-2\,^tGG_0\}.
\endalign$$
Here $(\det(-i\tau))^{1/2}$ is taken to be positive when $\tau=iY$, $Y>0$; in general,
the sign is found by analytic continuation.
\endproclaim

One can use this to derive the transformation formula (below), either as done in
Chapter 1 of [1], or alternatively, by adapting Eichler's argument (where $n$ was 1),
beginning with the identity
$$(A\tau+B)(C\tau+D)^{-1}=\,^tD^{-1}\,^tB+\,^tD^{-1}\tau(D\tau+D)^{-1}$$
(valid when $\det D\not=0$).  Eichler's approach yields the character as a generalised
Gauss sum, which can then be evaluated using fairly standard techniques (see [3]).  

\proclaim{Theorem 1.2 (Tranformation Formula)}  For $L$ a rank $m=2k+1$
lattice of level $N$ and 
$\gamma=\pmatrix
A&B\\C&D\endpmatrix\in\Gamma_0^{(n)}(N)$, 
$$\theta^{(n)}(L;\gamma\tau)
=\det(C\tau+D)^{m/2}\,\left(\frac{(-1)^k 2\det Q}{|\det D|}\right)
\,\sgn(\det D)^k
\,\theta^{(n)}(L;\tau).$$
Here, for $\det D\not=0$, $\lim_{\lambda\to0^+}\det(Ci\lambda I+D)^{1/2}=(\det D)^{1/2}$,
and in general the value of $\det(C\tau+D)^{1/2}$ is found by analytic continuation.
\endproclaim

In particular, this theorem applies to the basic Siegel theta series
$$\theta^{(n)}(\tau)=\theta^{(n)}((2);\tau)=\sum_{G\in\Z^{1,n}}\e\{2\,^tGG\tau\}$$
with $\gamma\in\Gamma_0^{(n)}(4)$.  Thus we can simplify our notation
by introducing automorphy factors and the slash operator, as follows.

Given $\gamma\in GSp_n^+(\Q)$, an automorphy factor for $\gamma$ is an
 analytic function $\varphi:\h_{(n)}\to \C$ so that 
$|\varphi(\gamma,\tau)|^2=|\det(C_{\gamma}\tau+D_{\gamma})|/\sqrt{\det\gamma}.$
Given $\gamma_1,\gamma_2\in GSp_n^+(\Q)$ with corresponding automorphy
 factors $\varphi_1,\varphi_2$, we define the product of the pairs
$[\gamma_1,\varphi_1(\tau)]$, $[\gamma_2,\varphi_2(\tau)]$ by
$$[\gamma_1,\varphi_1(\tau)]\,[\gamma_2,\varphi_2(\tau)]
=[\gamma_1\gamma_2,\varphi_1(\gamma_2\tau)\varphi_2(\tau)].$$
Correspondingly, we define a weight $m/2$ action of such a pair $[\gamma,\varphi(\tau)]$ on
$f:\h_{(n)}\to\C$ by
$$
f(\tau)|[\gamma,\varphi(\tau)]=f(\tau)|_{m/2}\,[\gamma,\varphi(\tau)]
=\varphi(\tau)^{-m}\,f(\gamma\tau).
$$
The default automorphy factor for $\gamma\in\Gamma_0^{(n)}(4)$ is
$${\theta^{(n)}(\gamma\tau)\over\theta^{(n)}(\tau)},$$
and we write $\widetilde\gamma$ to denote
$[\gamma,\theta^{(n)}(\gamma\tau)/\theta^{(n)}(\tau)].$

\smallskip

\noindent
{\bf Definition.}  Given $k, n, N\in\Z_+$ with $4|N$ and $\chi$ a Dirichlet
character modulo $N$, a Siegel modular form of degree $n$, weight $k+1/2$,
level $N$ and character $\chi$ is an analytic function $f:\h_{(n)}\to\C$
so that 
$$f|_{k+1/2}\,\widetilde\gamma=\chi(\det D_{\gamma})\,f$$
for all $\gamma\in\Gamma_0(N)$.  Here $D_{\gamma}$ denotes the lower right
$n\times n$ block of $\gamma$; often we simply write $\chi(\gamma)$ to
denote $\chi(D_{\gamma})$.
We write $\M_{k+1/2}(\Gamma_0^{(n)}(N),\chi)$ to denote the space of all such
modular forms.

\smallskip
\noindent{\bf Remark.}  Suppose $L,Q,N$ are as in Theorem 1.2.
Note that $\left(\theta^{(n)}(\tau)\right)^{m}=\theta^{(n)}(L_0;\tau)$
where $L_0$ is a rank $m$ lattice with quadratic form given by $2I_{m}$.
Thus for a lattice $L$ of odd rank $m$ and level $N$, and
$\gamma\in\Gamma_0^{(n)}(N)$, Theorem 1.2 gives us
$$\theta^{(n)}(L;\tau)|\widetilde\gamma=\chi(\det
D_{\gamma})\,\theta^{(n)}(L;\tau)$$
where, for $d\not=0$, $\chi(d)=\left(\frac{2\det Q}{|d|}\right).$

\smallskip

Suppose $f\in\M_{k+1/2}(\Gamma_0^{(n)}(N),\chi)$.  Then for any $B\in\Z^{n,n}_{\sym},$
we have $\pmatrix I&B\\&I\endpmatrix\in\Gamma_0^{(n)}(N)$, and hence
$f(\tau+B)=f(\tau)$.  Since $f$ is also analytic, $f$ has a Fourier series expansion
$$f(\tau)=\sum_T c(T) \e\{T\tau\}$$
where $T$ varies over $n\times n$, even integral matrices $T$ with
 $T\ge 0$ (meaning $T$ is positive semi-definite).  
Further, for any $G\in GL_n(\Z)$, we have $\gamma=\pmatrix G^{-1}\\&^tG\endpmatrix\in Sp_n(\Z)$
and $\theta^{(n)}(\gamma\tau)=\theta^{(n)}(\tau)$; so $f(G^{-1}\tau\,^tG^{-1})=
\chi(\det G) f(\tau),$ and consequently, $c(\,^tGTG)=\chi(\det G) c(T).$
Hence we can write 
$$f(\tau)=\sum_{\cls \Lambda}c(\Lambda) \e^*\{\Lambda\tau\}$$
where $\cls \Lambda$ varies over isometry classes of rank $n$ lattices equipped
with even integral, positive semi-definite quadratic form (oriented when $\chi(-1)=-1$),
$c(\Lambda)=c(T)$ where $T$ is a matrix for the quadratic form on $\Lambda$, and
$\e^*\{\Lambda\tau\}=\sum_G \e\{\,^tGTG\tau\}$ with $G\in O(T)\backslash GL_n(\Z)$
if $\chi(-1)=1$, $G\in O^+(T)\backslash SL_n(\Z)$ if $\chi(-1)=-1$.  (Here $O(T)$
is the orthogonal group of $T$, $O^+(T)=O(T)\cap SL_n(\Z)$.)

\noindent{\bf Definitions.}
A pair of $n\times n$ integer matrices $(C,D)$ is called a coprime symmetric
pair if $C\,^tD$ is symmetric, and for $G\in GL_n(\Q)$, $G(C,D)$ is integral
only if $G\in GL_n(\Z)$.  (Note that if 
$\pmatrix A&B\\C&D\endpmatrix\in Sp_n(\Z)$, then $(C,D)$ is a coprime
symmetric pair, as is $(\,^tB,\,^tD)$.  Conversely, if $(\,^tB,\,^tD)$ is a coprime
symmetric pair, then there is some $\pmatrix A&B\\C&D\endpmatrix\in Sp_n(\Z)$.)
For $(^tB,\,^tD)$ a coprime symmetric pair of $n\times n$ integral matrices with
$\det D\not=0$, 
 we define a generalised
Gauss sum
$$\G_B(D)=\sum_{G\in\Z^{1,n}/\Z^{1,n}\,^tD}\e\{2\,^tGGBD^{-1}\}.$$
Thus for $b,d$ coprime integers with $d\not=0$, $\G_b(d)$ is the classical
Gauss sum.
\smallskip

\noindent{\bf Remark.}  When following Eichler's approach to prove the
Transformation Formula, we encounter the Gauss sum 
$$\G_B(D;Q)=\sum_{G\in\Z^{m,n}/\Z^{m,n}\,^tD}\e\{\,^tGQGBD^{-1}\},$$
which we can evaluate using the theory of quadratic forms over finite
fields.

\proclaim{Proposition 1.3} Suppose $p$ is 
an odd prime, $D=\diag\{I_{r_0},pI_{r_1},p^2I_{r_2},I_{r_3}\}$
and 
$Y$ is the integral matrix
$$Y=\pmatrix Y_0&pY_2&0&Y_3\\^tY_2&Y_1&0\\0&0&I_{r_2}\\^tY_3\endpmatrix$$
where, for $i=0,1$, $Y_i$ is $r_i\times r_i$ and
symmetric.  Then
$$\G_Y(D)=p^{r_2}\,\left(\frac{\det Y_1}{p}\right)\,\G_1(p)^{r_1}$$
where $\G_1(p)$ is the classical Gauss sum.
\endproclaim

\demo{Proof}  
With $Y,D$ as above, we have
$$\G_Y(D)=\sum_{G_1\in\Z^{1,r_1}}\e\{2\,^tG_1G_1Y_1/p\}
\sum_{G_2\in\Z^{1,r_2}}\e\{2\,^tG_2G_2/p^2\}.$$
It is easy to see that
$$\sum_{G_2}\e\{2\,^tG_2G_2/p^2\}=\left(\sum_{g\in\Z/p\Z}\e\{2g^2/p^2\}\right)^{r_2}
=\G_1(p^2)^{r_2}=p^{r_2}.$$
Then, since $Y_1$ is symmetric and $p$ is odd, the quadratic form given by $Y_1$ can be
diagonalised over $\F_p$.  Since $SL_r(\Z)$ maps onto $SL_r(\F_p)$
(see, for instance, p. 21 of [12]), we can find
$E\in SL_{r_1}(\Z)$ so that $EY_1\,^tE\equiv\diag\{u_1,\ldots,u_{r_1}\}\ (\mod p)$, and thus
$$\sum_{G_1}\e\{2\,^tG_1G_1Y_1/p\}=\sum_{G_1}\e\{2\,^t(G_1E^{-1})(G_1E^{-1})U/p\}$$
where $U=\diag\{u_1,\ldots,u_{r_1}\}$.  Since $G_1E^{-1}$ varies over $\Z^{1,r_1}/p\Z^{1,r_1}$ as
$G_1$ does, we have
$$\sum_{G_1}\e\{2\,^tG_1G_1Y_1/p\}=\prod_{i=1}^{r_1}\G_{u_i}(p)=\left(\frac{u_1\cdots u_{r_1}}{p}\right)
(\G_1(p))^{r_1}.$$
Finally, note that $\det Y_1=\det U$, so $\left(\frac{u_1\cdots u_{r_1}}{p}\right)
=\left(\frac{\det Y_1}{p}\right)$.  $\square$
\enddemo

Next we introduce some notation for various types of representation numbers.
Let $p$ be a prime, $\F=\F_p$; let $T\in\F^{d,d}_{\sym}$, $S\in\F^{a,a}_{\sym}$, $a\le d$.
We set
$$\align
r(T,S)&=\#\{C\in\F^{d,a}:\ ^tCTC=S\ \},\\
r^*(T,S)&=\#\{C\in\F^{d,a}:\ ^tCTC=S,\ \rank_{\F}C=a\ \}.
\endalign$$
Thus $r(T,S)$ is the number of times $T$ represents $S$, and $r^*(T,S)$ is the
number of times $T$ primitively represents $S$.  We use $O(T)$ to denote the orthogonal
group of $T$, i.e.
$$O(T)=\{G\in GL_d(\F):\ ^tGTG=T\ \},$$ and set $o(T)=\# O(T)$.
When $Y=\ ^tGTG$ for some $G\in GL_d(\F)$, we write $Y\sim T$.

When $V$ is a dimension $d$ vector space over $\F$ with a quadratic form given by
a symmetric matrix $T$ (relative to some basis), we call $V$ a quadratic space over
$\F$, and we write $V\simeq T$.  We say $V$ is regular if $\det T\not=0$.
With $W$ another quadratic space over $\F$ with dimension $a\le r$ and quadratic form
given by a symmetric matrix $S$, we set
$$R^*(V,W)=\frac{r^*(T,S)}{o(S)};$$
so $R^*(V,W)$ is the number of dimension $a$ subspaces $W'$ of $V$ isometric to $W$,
meaning that relative to some basis for $W'$, the quadratic form on $W'$ is given
by $S$.  If $\dim V<\dim W$, then $R^*(V,W)=0$.

Now suppose the prime $p$ is odd.
We use $\H$ to denote the dimension 2 quadratic space over $\F_p$
with quadratic form given by the matrix $\big<1,-1\big>=\diag\{1,-1\}$; we also write
$\H\simeq\big<1,-1\big>$.  Similarly, $\A$ denotes the dimension 2 quadratic space over $\F_p$
so that $\A\simeq\big<1,-\omega\big>$, where $\left(\frac{\omega}{p}\right)=-1$.  
Given
any dimension $d$ quadratic space $V$ over $\F_p$, $V$ splits as $V_0\perp R$, where $V_0$
is regular and $R=\rad V\simeq\big<0\big>^s$, some $s\ge0$ (where $\big<0\big>^s$ denotes
the $s\times s$ matrix of zeros).  Also, the isometry class of $V_0$ is determined by the dimension
of $V_0$ and the value of the Legendre symbol $\left(\frac{\disc V_0}{p}\right)$ where
$\disc V_0$ is the determinant of a matrix for the quadratic form on $V_0$ (so $\disc V_0$ is
well-defined up to squares in $\F_p^{\times}$).  We have
$$V_0\simeq
\cases \H^c&\text{if $\dim V_0=2c$, $\left(\frac{\disc V_0}{p}\right)=\left(\frac{-1}{p}\right)^c$,}\\
\H^{c-1}\perp\A&\text{if $\dim V_0=2c$, $\left(\frac{\disc V_0}{p}\right)\not=\left(\frac{-1}{p}\right)^c$,}\\
\H^c\perp\big<1\big>
&\text{if $\dim V_0=2c+1$, $\left(\frac{\disc V_0}{p}\right)=\left(\frac{-1}{p}\right)^c$,}\\
\H^c\perp\big<\omega\big>
&\text{if $\dim V_0=2c+1$, $\left(\frac{\disc V_0}{p}\right)\not=\left(\frac{-1}{p}\right)^c$.}
\endcases$$
(Note that while $V_0$ is not uniquely determined by $V$, the isometry class of $V_0$ is.)
We say $V$ is totally isotropic if $V=\rad V$, i.e. $V\simeq\big<0\big>^s$.

We now define another generalised Gauss sum that we will later encounter.

\smallskip
\noindent{\bf Definition.} Suppose $p$ is an odd prime and $V$ is a dimension $d$ quadratic space
over $\F=\F_p$ with quadratic form given by $T$ modulo $p$, where $T$ is an even integral
matrix over $\Z$.  Then we define the twisted Gauss sum $\G^*(V)$ by
$$\G^*(V)=\sum_{Y\in\F^{d,d}_{\sym}} \left(\frac{\det Y}{p}\right)\e\{YT/p\}.$$
We define the normalised twisted Gauss sum $\widetilde \G(V)$ by
$$\widetilde\G(V)=p^{-d}(\G_1(p))^d \G^*(V).$$
When $V=\{0\}$, we agree that $\widetilde\G(V)=1$.
\smallskip

By Theorem 1.3 of [11], we have the following.

\proclaim{Proposition 1.4}  With $p$ an odd prime,  we have:
$$\align
\widetilde\G(\H^c\perp\big<0\big>^s)
&=\widetilde\G(\H^{c-1}\perp\A\perp\big<0\big>^s)\\
&={\cases (-1)^c p^{(c+x)^2-(c+x)}\prod_{i=1}^x (p^{2i-1}-1)&\text{if $s=2x$,}\\
0&\text{if $s=2x+1$;}\endcases}
\endalign$$
and with $\eta\in\F^{\times}$,
$$\widetilde\G(\H^c\perp\big<2\eta\big>\perp\big<0\big>^s)
=\cases (-1)^c\left(\frac{-\eta}{p}\right) p^{(c+x)^2+x}\prod_{i=1}^x (p^{2i-1}-1)
&\text{if $s=2x$,}\\
(-1)^c p^{(c+x)^2-(c+x)} \prod_{i=1}^x (p^{2i-1}-1)
&\text{if $s=2x-1$.}\endcases$$
\endproclaim

There are several elementary functions we use throughout:  For fixed prime $p$ and $m,r\in\Z$
with $r>0$,
$$\align
\delta(m,r)&=\delta_p(m,r)=\sum_{i=0}^{r-1}(p^{m-i}+1),\\
\mu(m,r)&=\mu_p(m,r)=\sum_{i=0}^{r-1}(p^{m-i}-1),\\
\delta(m,0)&=\mu(m,0)=1.
\endalign$$
When we have $m\ge r\ge 0$
$$\beta(m,r)=\beta_p(m,r)=\frac{\mu(m,r)}{\mu(r,r)},$$
which is the number of dimension $r$ subspaces of a dimension $m$ space over $\F_p$.
We will sometimes write, for instance, $\delta\mu(m,r)$ for $\delta(m,r)\cdot\mu(m,r)$.

We will frequently use the fact that $Tr(AB)=Tr(BA)$, and hence
$\e\{AB\}=\e\{BA\}$.

\bigskip
\head{\bf\S2. The standard generators of the Hecke algebra}\endhead
\smallskip

We begin this section
by defining the standard generators of the
Hecke algebra for half-integral weight Siegel modular
forms; then we analyse their action on Fourier coefficients.

Fix $N$ so that $4|N$, and set
$\widetilde\Gamma=\{\widetilde\gamma:\ \gamma\in\Gamma_0^{(n)}(N)\ \};$
let $\widetilde\delta=\left[\pmatrix pI\\&I\endpmatrix,p^{-n/2}\right]$.
Similar to the case of integral weight, we define
$$F|T(p)=\sum_{\widetilde\gamma}
\overline\chi(\gamma)\ F|\widetilde\delta^{-1}\widetilde\gamma$$
where $\widetilde\gamma$ runs over a complete set of representatives for
$(\widetilde\Gamma\cap\widetilde\delta\widetilde\Gamma\widetilde\delta^{-1})
\backslash\widetilde\Gamma$.
For $1\le j\le n$, set
$$X_j=\pmatrix pI_j\\&I_{n-j}\endpmatrix,\ \delta_j=\pmatrix X_j\\&X_j^{-1}\endpmatrix,
{\text{ and }}\widetilde\delta_j=[\delta_j,p^{-j/2}].$$
For $F\in\M_{k+1/2}(\Gamma^{(n)}_0(N),\chi)$, define
$$F|T_j(p^2)
=\sum_{\widetilde\gamma}\overline\chi(\gamma)\,F|\widetilde\delta_j^{-1}
\widetilde\gamma$$
where $\widetilde\gamma$ runs over a complete set of representatives for
$(\widetilde\Gamma\cap\widetilde\delta_j\widetilde\Gamma\widetilde\delta_j^{-1})
\backslash\widetilde\Gamma$.

\proclaim{Proposition 2.1} For $F\in\M_{m/2}(N,\chi)$ and $p$ prime, $F|T(p)=0$.
\endproclaim

\demo{Proof}  With $\ell\in2\Z_+$, write $\theta^{(n)}((\ell);\tau)$ for
$\theta^{(n)}(L;\tau)$ where $L=\Z x\simeq(\ell)$.
So $\theta^{(n)}(\tau)=\theta^{(n)}((2);\tau).$

Take $\gamma=\pmatrix A&B\\C&D\endpmatrix\in\Gamma=\Gamma_0^{(n)}(N)$ so that 
$\gamma'=\delta\gamma\delta^{-1}\in\Gamma$.
(So $p|C$, and $p\nmid\det D$.)  By Theorem 1.2,
$$\theta^{(n)}(\gamma\delta^{-1}\tau)
=\det(D\tau/p+D)^{1/2}\,\left(\frac{4}{|\det D|}\right)\,\theta^{(n)}(\tau/p),$$
and
$$\align
\theta^{(n)}(\gamma'\tau)&=\theta^{(n)}((2p);\gamma\tau/p)\\
&=\det(C\tau/p+D)^{1/2}\,\left(\frac{4p}{|\det D|}\right)\,\theta^{(n)}((2p);\tau/p)\\
&=\det(C\tau/p+D)^{1/2}\,\left(\frac{p}{|\det D|}\right)\,\theta^{(n)}(\tau).
\endalign$$
Thus 
$$\widetilde\delta \widetilde\gamma \widetilde\delta^{-1}
=\left[\gamma',\left(\frac{p}{|\det D|}\right) \theta^{(n)}(\gamma'\tau)/\theta^{(n)}(\tau)\right],$$
and hence with $\Gamma'=\delta\Gamma\delta^{-1}$,
$[\widetilde{(\Gamma\cap\Gamma')}:\widetilde\Gamma\cap\widetilde\delta
\widetilde\Gamma\widetilde\delta^{-1}]\le 2.$
To show this index is 2,
choose a prime $q\nmid N$ so that $\left(\frac{p}{q}\right)=-1$.
Choose $a,b\in\Z$ so that $aq-Nbp=1$, and set 
$$\gamma_0=\pmatrix a&&b\\&I_{n-1}&&0_{n-1}\\Np&&q\\&0_{n-1}&&I_{n-1}\endpmatrix\in\Gamma.$$
So 
$$\gamma'_0=\delta\gamma_0\delta^{-1}
=\pmatrix a&&pb\\&I_{n-1}&&0_{n-1}\\N&&q\\&0_{n-1}&&I_{n-1}\endpmatrix\in\Gamma,$$
but $\widetilde\gamma'_0\not=\widetilde\delta\widetilde\gamma_0\widetilde\delta^{-1}$.
Therefore $1,\widetilde\gamma_0'$ are coset representatives for
$(\widetilde\Gamma\cap\widetilde\delta\widetilde\Gamma\widetilde\delta^{-1})
\backslash(\widetilde{\Gamma\cap\Gamma'}),$ and
$$f|\widetilde\delta^{-1}\widetilde\gamma_0'
=-f|\widetilde\gamma_0\widetilde\delta^{-1}
=-\chi(\gamma_0)\,f|\widetilde\delta^{-1}.$$
So with $\widetilde\gamma$ running over a set of coset representatives for
$(\widetilde{\Gamma\cap\Gamma'})\backslash\widetilde\Gamma,$
and noting that $\chi(\gamma_0')=\chi(\gamma_0),$ we have
$$\align
f|T(p)
&=\sum_{\widetilde\gamma}\overline\chi(\gamma)\,f|\widetilde\delta^{-1}\widetilde\gamma
+\sum_{\widetilde\gamma}\overline\chi(\gamma_0'\gamma)\,
f|\widetilde\delta^{-1}\widetilde\gamma_0'\widetilde\gamma\\
&=\left(\sum_{\widetilde\gamma}\overline\chi(\gamma)\,f|\widetilde\delta^{-1}\widetilde\gamma\right)
- \left(\sum_{\widetilde\gamma}\overline\chi(\gamma)\,f|\widetilde\delta^{-1}\widetilde\gamma\right)\\
&=0,
\endalign$$
proving the proposition. $\square$
\enddemo

We use a similar argument to prove the following.

\proclaim{Lemma 2.2}  Let $\Gamma=\Gamma_0(N)$ and let
$\delta_j, \widetilde\delta_j$ be as above.  Set
$\Gamma_j'=\delta_j\Gamma\delta_j^{-1}$; then for $p\nmid N$,
$$\widetilde{\Gamma'_j\cap\Gamma}=
\widetilde\delta_j\widetilde\Gamma\widetilde\delta_j^{-1}\cap\widetilde\Gamma.$$
\endproclaim

\demo{Proof} 
Say $\gamma=\pmatrix A&B\\C&D\endpmatrix\in\Gamma$ so that
$\gamma'=\delta_j\gamma\delta_j^{-1}\in\Gamma$.  
Then by Theorem 1.2,
$$\align
\theta^{(n)}(\gamma'\tau)/\theta^{(n)}(\tau)&=\det(X_j^{-1}CX_j^{-1}\tau+X_j^{-1}DX_j)^{1/2}\\
&=\det(CX_j^{-1}\tau X_j^{-1}+D)^{1/2}\\
&=\theta^{(n)}(\gamma\delta_j^{-1}\tau)/\theta^{(n)}(\delta_j^{-1}\tau),
\endalign$$
and hence
$\widetilde\gamma'=\widetilde\delta_j\widetilde\gamma\widetilde\delta_j^{-1}.$
$\square$
\enddemo

\proclaim {Theorem 2.3}  Take $p$ a prime and
$F\in\M_{k+1/2}\left(\Gamma_0^{(n)}(N),\chi\right).$
\itemitem{(a)} If $p\nmid N$, then
$$F|T_j(p^2)=\sum_{\Omega,\Lambda'_1,Y}
\chi(\det D)
F|\widetilde\delta_j^{-1}
\left[\pmatrix D&^tY\\&D^{-1}\endpmatrix,{\G_Y(D)\over (\det D)}\right]
\widetilde{\pmatrix G^{-1}\\&^tG\endpmatrix}$$
where $\Omega, \Lambda'_1$ vary subject to
$p\Lambda\subseteq\Omega\subseteq{1\over p}\Lambda$, and
$\overline\Lambda'_1$ varies over all codimension $n-j$ subspaces of
$\Lambda\cap\Omega/p(\Lambda+\Omega)$. Here 
$$D=D(\Omega)=\diag\{I_{n_0},pI_{j-r},p^2I_{n_2},I_{n-j}\}$$ 
with $r=n_0+n_2$, and $G=G(\Omega,\Lambda'_1)\in
SL_n(\Z)$ 
so that
$$\Omega=\Lambda GD^{-1}X_j,\ 
\Lambda'_1=\Lambda G\pmatrix 0_{r_0}\\&I_{j-r}\\&&0\endpmatrix;
\ ^tY=\pmatrix Y_0&Y_2&0&Y_3\\p\
^tY_2&Y_1&0\\0&0&I\\^tY_3\endpmatrix$$
with $Y_0$ symmetric, $n_0\times n_0$, varying modulo $p^2$, $Y_1$
symmetric, $(n-r)\times (n-r)$ varying modulo $p$ with the restriction that $p\nmid\det
Y_1$, $Y_2$ $n_0\times (j-r)$, varying modulo $p$, $Y_3$ $n_0\times(n-j)$,
varying modulo $p$.
\itemitem{(b)} If $p|N$, then 
$$F|T_j(p^2)=\sum_{\Omega,Y}
F|\widetilde\delta_j^{-1}
\widetilde{\pmatrix I&Y\\&I\endpmatrix}
\widetilde{\pmatrix G^{-1}\\&^tG\endpmatrix}$$
where $\Omega$ varies subject to
$p\Lambda\subseteq\Omega\subseteq\Lambda$, $[\Lambda:\Omega]=p^j$,
$G=G(\Omega)\in
SL_n(\Z)$  so that
$\Omega=\Lambda GX_j$
and $$Y=\pmatrix Y_0&Y_3\\^tY_3&0\endpmatrix$$
with $Y_0$ symmetric, $j\times j$, varying modulo $p^2$, $Y_3$
$j\times(n-j)$, varying modulo $p$.
\endproclaim

\demo{Proof}  
By Lemma 2.2, a set of coset representatives for
the action of the half-integral weight Hecke operator
$T_j(p^2)$ is $\{\widetilde\gamma\}$ where $\{\gamma\}$ is a
set of coset representatives for the integral weight Hecke operator
$T_j(p^2)$, and a set of representatives for this was given in
\S6 of [6]; note that the matrices
$G$ presented there can be chosen from $SL_n(\Z)$.
Thus a set of coset representatives corresponding to $T_j(p^2)$ on
$\M_{k+1/2}(\Gamma^{(n)}_0(N),\chi)$ is
$$\left\{
\widetilde{\pmatrix D&^tY\\U&W\endpmatrix}
\widetilde{ \pmatrix G^{-1}\\&^tG\endpmatrix}
\right\}$$
where 
$\pmatrix D&^tY\\U&W\endpmatrix\in\Gamma=\Gamma^{(n)}_0(N)$ 
with $U=\diag\{U_0,0_{n-j}\}$, $W=\diag\{W_0,I_{n-j}\}$, and
$G=G(\Omega,\Lambda_1'),\ D=D(\Omega),\ Y$ vary as claimed.
(From the definition, we know that for $\gamma_1,\gamma_2\in\Gamma$, 
$\widetilde{\gamma_1\gamma_2}
=\widetilde\gamma_1\widetilde\gamma_2.$) 
Note that $DW\equiv I\ (\mod N)$, so $\overline\chi(\det W)
=\chi(\det D)$.

Suppose first $p\nmid N$. Set $X=X_j$, $X'=XUD^{-1}X$,
$$\gamma=\pmatrix D&^tY\\U&W\endpmatrix,\ 
\beta=\pmatrix D&^tY\\&D^{-1}\endpmatrix,\ 
\gamma'=\delta_j^{-1}\gamma\beta^{-1}\delta_j 
=\pmatrix I\\X'&I\endpmatrix;$$
since $XU=UX$, $X'$ is integral and divisible by $N$, and so $\gamma'\in\Gamma^{(n)}_0(N)$.
We will show that
$\widetilde\gamma'
=\widetilde\delta_j^{-1}\widetilde\gamma[\beta^{-1},(\det D)(\G_Y(D))^{-1}]
\widetilde\delta_j,$
and hence
$$
F|\widetilde\delta_j^{-1}\widetilde\gamma
=F|\widetilde\gamma'\widetilde\delta_j^{-1}[\beta,(\det D)^{-1}\G_Y(D)]
=F|\widetilde\delta_j^{-1}[\beta,(\det D)^{-1}\G_Y(D)].$$

We have
$$\widetilde\delta_j^{-1}\widetilde\gamma
\left[\beta^{-1},\frac{\det D}{\G_Y(D)}\right]
\widetilde\delta_j=
\left[\gamma',\frac{\det D}{\G_Y(D)}\cdot
{\theta^{(n)}(\gamma\beta^{-1}\delta_j\tau)\over\theta^{(n)}(\beta^{-1}\delta_j\tau)}\right],\ 
\beta^{-1}=\pmatrix D^{-1}&-Y\\&D\endpmatrix.$$
Also, since $DU$ must be symmetric, $DUD^{-1}$ is integral and divisible by 4.  Therefore,
using the Inverstion Formula (Theorem 1.1), we have
$$\align
&\theta^{(n)}(\gamma\beta^{-1}\delta_j\tau)\\
&\quad=
\frac{1}{\sqrt2^n}(\det X)^{-1}
 \left(\det(-i\tau(X'\tau+I)^{-1})\right)^{-1/2}\\
&\qquad\cdot
\sum_{g\in\Z^{1,n}}\e\left\{-\frac{1}{ 2}\,^tgg(UD^{-1}X\tau+X^{-1})
\tau^{-1}X^{-1}\right\}\\
&\quad=
\frac{1}{\sqrt2^n} 
(\det X)^{-1}
 \left(\det(-i\tau(X'\tau+I)^{-1})\right)^{-1/2}\\
&\qquad\cdot
\sum_{g_0\, (D)}\e\left\{-\frac{1}{2}\,^tg_0g_0UD^{-1}\right\}
\theta^{(n)}\left(\left(\frac{1}{2}\right),g_0D^{-1};-DX^{-1}\tau^{-1}X^{-1}D\right)\\
&\quad=
(\det X)^{-1} (\det(-i D^{-1}X\tau XD^{-1}))^{1/2}
\left(\det(-i\tau(XUD^{-1}X\tau+I)^{-1})\right)^{-1/2}\\
&\qquad\cdot
\sum_{g_0\,(D)}\e\left\{-{1\over 2}\,^tg_0g_0UD^{-1}\right\}
\sum_{g\in\Z^{1,n}}\e\{2\,^tggD^{-1}X\tau XD^{-1}-2D^{-1}\,^tg_0g\}.
\endalign$$
(By $g_0\, (D)$ we really mean $g\in\Z^{1,n}/\Z^{1,n}\ D$.)
Since $\gamma,\gamma^{-1}$ are symplectic, we know $UY=WD-I$,
$YU=\,^tWD-I$, and $D^{-1}\,^tU=UD^{-1}$ with $4|U$; thus
$$\align
& \e\{2\,^t(g_0U/2+g)(g_0U/2+g)YD^{-1}\}\\
&\quad = \e\left\{-{1\over 2}\,^tg_0g_0UD^{-1}-2D^{-1}\,^tg_0g\}\ \e\{2\
^tggYD^{-1}\right\}.
\endalign$$
Thus, remembering that $\det D,\det X >0$,
$$\align
\theta^{(n)}(\gamma\beta^{-1}\delta_j\tau)
&=
(\det D)^{-1} (\det(-i \tau))^{1/2}
\left(\det(-i\tau(XUD^{-1}X\tau+I)^{-1})\right)^{-1/2}\\
&\quad\cdot
\sum_{{g,g_0\in\Z^{1,n}}\atop{g_0\,(D)}}
\e\{2\,^t(g_0U/2+g)(g_0U/2+g)YD^{-1}\}\\
&\quad\cdot
\e\{2\,^tgg(D^{-1}X\tau XD^{-1}-YD^{-1})\}.
\endalign$$
We know that $(D,\,^tU)$ is a symmetric pair with
$U=\diag\{U_0,0_{n-j}\},$
so 
$$U=\pmatrix
W_0&W_3&W_5\\p\,^tW_3&W_1&W_4\\p^2\,^tW_5&p\,^tW_4&W_2\endpmatrix$$
with the diagonal blocks symmetric and of sizes $n_0\times n_0$,
$(j-r)\times(j-r)$, and $n_2\times n_2$.  Since we also have that
$(D,\,^tU)$ is a coprime pair, $p\nmid\det\pmatrix
W_1&W_4\\p\,^tW_4&W_2\endpmatrix;$
thus for fixed $g\in\Z^{1,n}$, $g_0U/2+g$ varies over $\Z^{1,n}/\Z^{1,n}D$
as $g_0$ does (recall that $4|U$).  Thus
the sum on $g_0$ is independent of the choice of $g$, and so the sum
on $g_0$ is $\G_Y(D)$. 
Hence 
$$\align
\theta^{(n)}(\gamma\beta^{-1}\delta_j\tau)
&=
(\det D)^{-1} (\det(-i \tau))^{1/2}
\left(\det(-i\tau(XUD^{-1}X\tau+I)^{-1})\right)^{-1/2}\\
&\quad\cdot \G_Y(D)\,\theta^{(n)}(\beta^{-1}\delta_j\tau).
\endalign$$

Next, recall that $XU=UX$ with $N|X'$, and hence
for $g\in\Z^{1,n}$, $4|\,^tggX'$.  Thus,
using Theorem 1.1, we have
$$\align
\theta^{(n)}(\gamma'\tau)
&=\frac{1}{\sqrt 2^n} (\det(-i\tau(X'\tau+I)^{-1}))^{-1/2}
\sum_{g\in\Z^{1,n}}\e\left\{-\frac{1}{2}\,^tgg\tau^{-1}\right\}\\
&=(\det(-i\tau(X'\tau+I)^{-1}))^{-1/2} (\det(-i\tau))^{1/2}
\theta^{(n)}(\tau).
\endalign$$
This means
$$
\frac{\theta^{(n)}(\gamma'\tau)}{\theta^{(n)}(\tau)}
=\frac{\G_Y(D)}{\det D}
\frac{\theta^{(n)}(\gamma\beta^{-1}\delta_j\tau)}{\theta^{(n)}(\beta^{-1}\delta_j\tau)},
$$
completing the proof in the case that $p\nmid N$.

In the case $p|N$, the argument is much simpler, as
the coset representatives for
$(\Gamma'_j\cap\Gamma)\backslash\Gamma$ are those representatives as
above where $D=I$.
Since $\Omega=\Lambda GX$, we have
$p\Lambda\subseteq\Omega\subseteq\Lambda$ with $[\Lambda:\Omega]=p^j$.
$\square$
\enddemo

We complete this section by evaluating the action of $T_j(p^2)$ on the Fourier
coefficients of a half-integral weight Siegel modular form.  These involve the
normalised twisted Gauss sums, as defined in \S1 and whose values are given in
Proposition 1.4.

\proclaim{Theorem 2.4}  Take $f\in\M_{k+1/2}(\Gamma^{(n)}_0(N),\chi)$ where
$4|N$, and let $p$ be a prime.  Let $\chi'$ be the character modulo $N$ defined by
$$\chi'(d)=\chi(d)\left({(-1)^{k+1}\over |d|}\right)(\sgn\,d)^{k+1}.$$
\itemitem{(a)}  Suppose $p\nmid N$.  
Given
$$\Lambda=\Lambda_0\oplus\Lambda_1\oplus\Lambda_2,\ 
\Omega=p\Lambda_0\oplus\Lambda_1\oplus {1\over p}\Lambda_2$$
with $n_i=\rank\Lambda_i$, $r=n_0+n_2$, set
$$A_j(\Lambda,\Omega)=\chi'(p^{j-r})p^{j/2+k(n_2-n_0)+n_0(n-n_2)}
\sum_{{\cls U}\atop{\dim
U=j-r}}R^*(\Lambda_1/p\Lambda_1,U)\widetilde\G(U)$$
if $\Lambda, \Omega$ are even integral, and set $A_j(\Lambda,\Omega)=0$
otherwise.
Then the $\Lambda$th coefficient of
$f|T_j(p^2)$ is
$$\sum_{p\Lambda\subseteq\Omega\subseteq{1\over p}\Lambda}
A_j(\Lambda,\Omega) c(\Omega).$$
\itemitem{(b)} Suppose $p|N$.  Then the $\Lambda$th coefficient of $f|T_j(p^2)$ is
$$p^{j(n-k+1/2)}\sum_{{p\Lambda\subseteq\Omega\subseteq\Lambda}\atop
{[\Lambda:\Omega]=p^j}}
c(\Omega).$$
\endproclaim

\demo{Proof}  
Suppose first that
$p\nmid N$; then $p\not=2$ and 
$$\align
f(\tau)|T_j(p^2)&=p^{-j(k+1/2)}\sum_{D,Y,G}\chi(\det D)(\det
D)^{2k+1} \G_Y(D)^{-2k-1}\\
&\quad \cdot
\sum_T
c(T)\e\{TX^{-1}DG^{-1}\tau\,^tG^{-1}DX^{-1}\}\,\e\{TX^{-1}\,^tYDX^{-1}\}
\endalign$$
where $D,Y,G$ vary as in Theorem 2.3, and $T$ varies over all $n\times
n$ even integral, positive semi-definite matrices.

Fix $T,G$ and $D=\diag\{I_{n_0},pI_{j-r},p^2I_{n_2},I_{n-j}\},$ 
and let $Y$ vary.
As described in Theorem 2.3, we have
$$^tY=\pmatrix Y_1&Y_2&0&Y_3\\p\,^tY_2&Y_1&0\\0&0&I\\^tY_3\endpmatrix;
\text{ write }
T=\pmatrix
T_0&T_2&*&T_3\\^tT_2&T_1&*&*\\*&*&*&*\\^tT_3&*&*&*\endpmatrix$$
with $T_i$ the size of $Y_i$.
By Proposition 1.3, $\G_Y(D)=p^{n_2}\left(\frac{\det Y_1}{p}\right)\G_1(p)^{j-r}$; so
$$\align
&\sum_Y \G_Y(D)^{-2k-1}\e\{TX^{-1}\,^tYDX^{-1}\}\\
&\quad
= p^{-n_2(2k+1)} \G_1(p)^{(r-j)(2k+1)} \sum_{Y_0\,(p^2)}\e\{T_0Y_0/p^2\}
\cdot\sum_{Y_2\,(p)}\e\{2T_2Y_2/p\}\\
&\qquad\cdot
\sum_{Y_3\,(p)}\e\{2T_3Y_3/p\}\cdot 
\sum_{Y_1\,(p)}\left(\frac{\det Y_1}{p}\right)\e\{T_1Y_1/p\}.
\endalign$$
If $T_0\equiv0\,(\mod p^2)$, $T_2\equiv0\,(\mod p)$, 
$T_3\equiv0\,(\mod p)$ then
the sum on $Y$ is
$$\left(\frac{-1}{p}\right)^{(j-r)(k+1)} p^{-n_2(2k+1)+(r-j)k+n_0(n-n_2+1)}
\widetilde\G(T_1\ \mod p),$$
and otherwise the sum on
$Y$ is 0  (here we used that $(\G_1(p))^2=\left(\frac{-1}{p}\right)p$.)
Therefore,
$$\align
f(\tau)|T_j(p^2)
&=
\sum_{D,G}\chi'(\det D)
p^{j/2+k(n_2-n_0)+n_0(n-n_2)}\\
&\qquad\cdot
 \sum_{{T}\atop{T[X^{-1}D]\text{ integral}}} 
\widetilde\G(T_1\ \mod p) c(T) \e\{T[X^{-1}DG^{-1}]\tau\}.
\endalign$$
Let $\Lambda=\Z x_1\oplus\cdots\oplus\Z x_n$ be equipped with the quadratic
form $T[X^{-1}DG^{-1}]$ (relative to the given basis); then
$\Omega=\Lambda GD^{-1}X\simeq T$ relative to 
$(x_1\,\ldots\,x_n)GD^{-1}X$.  Also, relative to these bases we have splittings
$$\Lambda=\Lambda_0\oplus\Lambda'_1\oplus\Lambda_2\oplus\Lambda_1'',\ 
\Omega=p\Lambda_0\oplus\Lambda_1'\oplus{1\over p}\Lambda_2\oplus\Lambda_1''$$
where $n_0=\rank\Lambda_0$, $j-r=\rank\Lambda'_1$, $n_2=\rank\Lambda_2$, and
$\Lambda_1'\simeq T_1$.  
From [HW] we know that with $\Lambda_1=\Lambda_1'\oplus\Lambda_1''$
and $\Omega$ fixed, $G=G(\Omega,\Lambda_1')$ varies to vary
$\overline\Lambda_1'$ over all dimension $j-r$ subspaces of
$\Lambda_1/p\Lambda_1\approx \Lambda\cap\Omega/p(\Lambda+\Omega)$.
Thus
$$\sum_{\Lambda'_1}\widetilde\G(T_1\ \mod p)
=\sum_{{\cls U}\atop{\dim U=j-r}}R^*(\Lambda_1/p\Lambda_1,U)\widetilde\G(U).$$
This proves the theorem in the
case $p\nmid N$.

In the case $p|N$, the analysis is simpler, as $D$ is always $I$ (and
hence $n_0=r=j$, $n_2=0$);
following the above reasoning yields the theorem in this case.
$\square$
\enddemo

\bigskip
\head{\bf \S3. An alternate set of generators for the Hecke algebra}\endhead
\smallskip

In pursuit of a more utile formula for the action of Hecke operators
on Fourier coefficients, we introduce an alternate set of generators.
When we studied the action of Hecke operators on integral
weight Siegel modular
forms in [6], in the case $p\nmid N$ we encountered incomplete character sums,
which we completed these sums by replacing $T_j(p^2)$ with
$$\widetilde T_j(p^2)=p^{j(k-n-1)}\sum_{\ell=0}^j\chi(p^{j-\ell})
\,\beta(n-\ell,j-\ell)\,T_{\ell}(p^2).$$
Here, with half-integral weight, the incomplete character sums are twisted,
giving us the generalised twisted Gauss sums.
In Theorem 3.3 we show that by replacing $T_j(p^2)$ with
$$\widetilde T_j(p^2)=p^{j(k-n)}\sum_{\ell=0}^j p^{-\ell/2}\,\chi'(p^{j-\ell})
\,\beta(n-\ell,j-\ell)\,T_{\ell}(p^2),$$
we eliminate these Gauss sums from the formula for the action on Fourier coefficients.

To ready ourselves to prove Theorem 3.3, we first establish the following
relationship between $\widetilde \G(W)$ and a weighted sum of representation numbers.

\proclaim{Lemma 3.1}  For $p$ an odd prime, $\F=\F_p$, and
$W$ a quadratic space over $\F$ with dimension
$m\ge 0$,
$$\widetilde \G(W)=\sum_{a=0}^m (-1)^{m+a} p^{m(m-1)/2+a(a-m)}
R^*(W\perp\big<2\big>,\big<0\big>^a).$$
\endproclaim

\demo{Proof}  
Since $R^*(\big<2\big>,\{0\})=1$, the lemma holds for $m=0$.
So suppose $m\ge 1$.  Then standard theory tells us that with $\F v\simeq\big<2\big>$,
$W\perp\F v$ splits as $W_0\perp R$ where $R=\rad (W\perp\F v)\simeq\big<0\big>^s$,
some $s\in\Z_{\ge 0}$, with $W_0$ regular.  Any totally isotropic subspace $U$ of $W\perp\F v$ splits as
$U_0\perp U_1$ where $U_1=U\cap R$.  Given a dimension $t$
subspace $U_1$ of $R$
and $a\ge t$, the number of distinct totally isotropic, dimension $a$ subspaces
$U$ of $W\perp\F v$ with $U\cap R=U_1$ is
$$p^{(s-t)(a-t)} R^*(W_0\perp\big<2\big>,\big<0\big>^{a-t}).$$
Since $\beta(s,t)$ is the number of dimension $t$ subspaces of $R$,
$$R^*(W\perp\big<2\big>,\big<0\big>^a)
=\sum_{t=0}^s \beta(s,t)\, p^{(s-t)(a-t)}\,R^*(W_0\perp\big<2\big>,\big<0\big>^{a-t}).$$

For $s,t,c,q\in\Z$ with $s,t,c\ge 0$, set
$$\align
S_t(c,q)
&=(-1)^c p^{c(c-t+q)} \sum_{\ell=0}^c
(-1)^{\ell} p^{\ell(\ell-2c+t-q)} \beta(c,\ell)\delta(c+q,\ell),\\
X_s(c,q)&=\sum_{t=0}^s(-1)^tp^{t(t-c-s-q)}\beta(s,t)S_t(c,q).
\endalign$$
Using, for instance, Theorems 2.59 and 2.60 of [5], we know
$$R^*(W_0\perp\big<2\big>,\big<0\big>^{\ell})
=\cases
\beta(c,\ell)\delta(c-1,\ell)&\text{if $W_0\perp\big<2\big>\simeq\H^c,$}\\
\beta(c-1,\ell)\delta(c,\ell)&\text{if $W_0\perp\big<2\big>\simeq\H^{c-1}\perp\A,$}\\
\beta(c,\ell)\delta(c,\ell)&\text{if $W_0\perp\big<2\big>\simeq\H^{c}\perp\big<\eta\big>.$}
\endcases$$
(Here $\eta\in\F^{\times}$.)  Hence, with $\omega\in\F^{\times}$ so that
$\left(\frac{\omega}{p}\right)=-1$, 
$$\align
&\sum_{a=0}^m (-1)^a\,p^{a(a-m)}\,R^*(W\perp\big<2\big>,\big<0\big>^a)\\
&\quad=
{\cases
(-1)^c\,p^{c(1-c)}\,X_s(c,-1)&\text{if $W\simeq\H^{c-1}\perp\big<-2\big>\perp\big<0\big>^s,$}\\
(-1)^{c-1}\,p^{c(1-c)}\,X_s(c-1,1)&\text{if $W\simeq\H^{c-1}\perp\big<-2\omega\big>\perp\big<0\big>^s,$}\\
(-1)^c\,p^{-c^2}\,X_s(c,0)&\text{if $W\simeq W_0\perp\big<0\big>^s,\ \dim W_0=2c.$}
\endcases}
\endalign$$

When $c>0$, replacing $\ell$ by $c-\ell$ and using the identities $\beta(c,\ell)=\beta(c,c-\ell)$
and $\beta(c,\ell)=p^{\ell}\beta(c-1,\ell)+\beta(c-1,\ell-1)$,
and then replacing $\ell-1$ by $\ell$ we get
$$S_0(c,q)=S_0(c-1,q+1)=S_0(0,q+c)=1.$$
Clearly $S_t(0,q)=1$; using the definitions of $\beta$ and $\delta$, when $c>0$ we have
$$\align
&S_t(c,q)+(p^c-1)(p^{c+q}+1)\,S_t(c-1,q)\\
&\quad=(-1)^c p^{c(c-t+q)}\sum_{\ell=0}^c (-1)^{\ell}p^{\ell(\ell-2c+t-q)}
\frac{\mu(c,\ell)\delta(c+q,\ell)}{\mu(\ell,\ell)}\\
&\qquad
+(-1)^c p^{c(c-t+q)-2c-q+t+1}\sum_{\ell=1}^c (-1)^{\ell} p^{(\ell-1)(\ell+1-2c+t-q)}
\frac{\mu(c,\ell)\delta(c+q,\ell)}{\mu(\ell-1,\ell-1)}\\
&\quad=(-1)^c p^{c(c-t+q)}\left[ 1+\sum_{\ell=1}^c (-1)^{\ell}p^{\ell(\ell+1-2c+t-q)}
\beta(c,\ell)\delta(c+q,\ell)\right] \\
&\quad=p^cS_{t+1}(c,q).
\endalign$$
Taking $S_t(c,q)=0$ when $c<0$, the above relation also holds for $c=0$.

Using that $\beta(s+1,t)=p^t\beta(s,t)+\beta(s,t-1),$ and the recursion relation
for $S_t(c,q)$, we have
$$\align
X_{s+1}(c,q)
&=\sum_{t=0}^s (-1)^t p^{t(t-c-s-q)} \beta(s,t) S_t(c,q)\\
&\quad + \sum_{t=1}^{s+1} (-1)^t p^{t(t-c-s-q-1)} \beta(s,t-1) S_t(c,q)\\
&=X_s(c,q)-p^{-2c-s-q}X'_s(c,q)-p^{-2c-s-q}(p^c-1)(p^{c+q}+1)X_s(c-1,q)
\endalign$$
where
$$X'_s(c,q)=\sum_{t=0}^s(-1)^t p^{t(t-c-s-q+1)}\,\beta(s,t)\,S_t(c,q).$$
Similarly, using that $\beta(s+1,t)=\beta(s,t)+p^{s+1-t}\beta(s,t-1)$,
$$X'_{s+1}(c,q)=X_s(c,q)-p^{-2c+1-q}X'_s(c,q)-p^{-2c+1-q}(p^c-1)(p^{c+q}+1)X_s(c-1,q).$$

Using induction on $x$,
the recursion relations for $X$ and $X'$, and the fact that $S_0(c,q)=1$, we get
the following.
$$\align
X_s(c,-1)
&={\cases p^{2x-2cs-x^2}\prod_{i=1}^x(p^{2i-1}-1)&\text{if $s=2x$,}\\
-p^{-c-2cx-x^2}\prod_{i=1}^{x+1}(p^{2i-1}-1)&\text{if $s=2x+1$,}
\endcases}\\
X'_s(c,-1)
&={\cases 1&\text{if $s=0$,}\\
X_{s-1}(c,-1)+(p^{c-1}+1)X_s(c-1,-1)&\text{if $s=2x>0$,}\\
(p^{c-1}+1)X_s(c-1,-1)&\text{if $s=2x+1$;}
\endcases}\\
X_s(c-1,1)
&={\cases p^{2x-2cs-x^2}\prod_{i=1}^x(p^{2i-1}-1)&\text{if $s=2x$,}\\
p^{-c-2cx-x^2}\prod_{i=1}^{x+1}(p^{2i-1}-1)&\text{if $s=2x+1$,}
\endcases}\\
X'_s(c-1,1)
&={\cases 1&\text{if $s=0$,}\\
X_{s-1}(c-1,1)+(p^{c-1}-1)X_s(c-2,1)&\text{if $s=2x>0$,}\\
-(p^{c-1}-1)X_s(c-2,1)&\text{if $s=2x+1$;}
\endcases}\\
X_s(c,0)&={\cases p^{-2cx-x^2}\prod_{i=1}^x(p^{2i-1}-1)&\text{if
$s=2x$,}\\
0&\text{if $s=2x+1$;}
\endcases}\\
X'_s(c,0)&{=\cases 1&\text{if $s=0$,}\\
X_s(c-1,0)&\text{if $s=2x>0$,}\\
-p^{-2cx-x^2}\prod_{i=1}^{x+1}(p^{2i-1}-1)&\text{if $s=2x+1$.}
\endcases}
\endalign$$
The lemma now follows from these identities and Proposition 1.4. $\square$
\enddemo

Using this, we also establish the following.

\proclaim{Lemma 3.2} Suppose $W$ is a dimension $m\ge0$ quadratic space over
$\F=\F_p$, $p$ an odd prime.
\item{(a)}  For $0\le a\le m$,
$$R^*(W\perp\big<2\big>,\big<0\big>^a)
=R^*(W,\big<0\big>^a)+2R^*(W,\big<0\big>^{a-1}\perp\big<-2\big>).$$
\item{(b)}  With $\cls U$ denoting the isometry class of a space $U$,
$$\sum_{q=0}^m \sum_{{\cls U}\atop{\dim U=q}}
R^*(W,U)\,\widetilde\G(U)
=p^{m(m-1)/2} R^*(W\perp\big<2\big>,\big<0\big>^m).$$
\endproclaim

\demo{Proof} (a) Fix $a$, $0\le a\le m$.  If $a=0$, then the claim is that
$$R^*(W\perp\big<2\big>,\{0\})=R^*(W,\{0\}),$$
which holds since there is only 1 dimension 0 subspace of any space.
So suppose $1\le a\le m$.  Let $U'$ be a dimension $a$ totally isotropic
subspace of $W\perp\F v$ where $\F v\simeq\big<2\big>$, and let $U$
be the projection of $U'$ onto $W$.  Then either $U=U'\simeq\big<0\big>^a$,
or $U\simeq\big<0\big>^{a-1}\perp\big<-2\big>$; also, there are exactly
2 totally isotropic subspaces of $W\perp\F v$ that project onto a
given subspace $U\simeq\big<0\big>^{a-1}\perp\big<-2\big>$ of $W$, and from this
the claim in (a) follows.

(b)  Here we argue by induction on $m$.  For $m=0$, the claim is trivially true.
So suppose $m\ge 1$.  Using the identity
$\beta(m,a)=\beta(m-1,a)+p^{m-a}\beta(m-1,a-1),$ and then replacing $a-1$ by $a$,
we get
$$\sum_{a=0}^m (-1)^{m+a}p^{(m-a)(m-a-1)/2}\beta(m,a)=0.$$
Thus for $a\le q\le m$, replacing $m$ and $a$ by $m-q$ and $a-q$ in the above
identity, we get 
$$\sum_{a=q}^{m-1}(-1)^{m+a+1}p^{(m-a)(m-a-1)/2}\beta(m-q,a-q)=1.$$

With $W$ a dimension $m$ space over $\F$,
the number of ways to extend a dimension $q$ subspace $U$ of 
 $W$ to a dimension $a$ subspace $Y$ of $W$ is $\beta(m-q,a-q)$.  Thus
with the preceding expression for 1,  using the induction hypothesis
and (a) of this lemma, we have
$$\align
&\sum_{{\cls U}\atop{\dim U<m}} R^*(W,U)\widetilde\G(U)\\
&\qquad =
\sum_{a=0}^{m-1}(-1)^{m+a+1}p^{(m-a)(m-a-1)/2}
\sum_{{\cls Y}\atop{\dim Y=a}}R^*(W,Y)
\sum_{{\cls U}\atop{\dim U\le a}}R^*(Y,U)\widetilde\G(U)\\
&\qquad =
\sum_{a=0}^{m-1}(-1)^{m+a+1}p^{m(m-1)/2+a(a-m)}R^*(W\perp\big<2\big>,\big<0\big>^a).
\endalign$$
Now add $\widetilde\G(W)=R^*(W,W)\widetilde\G(W)$ to both sides of
this equation; using Lemma 3.1 yields the result.
$\square$
\enddemo

Now it is easy to prove our main formula.

\proclaim{Theorem 3.3} Take $f\in\M_{k+1/2}(\Gamma_0^{(n)}(N),\chi)$ where
$4|N$, and let $p$ be a prime such that $p\nmid N$; let $\chi'$ be defined as in 
Theorem 2.4.
Given
$$\Lambda=\Lambda_0\oplus\Lambda_1\oplus\Lambda_2,\ 
\Omega=p\Lambda_0\oplus\Lambda_1\oplus {1\over p}\Lambda_2$$
with $n_i=\rank\Lambda_i$, $r=n_0+n_2$, set
$$E_j(\Lambda,\Omega)=j(k-n)+k(n_2-n_0)+n_0(n-n_2)+(j-r)(j-r-1)/2;$$
set
$$\widetilde A_j(\Lambda,\Omega)
=\chi'(p^{j-r})p^{E_j(\Lambda,\Omega)}
R^*(\Lambda_1/p\Lambda_1\perp\big<2\big>,\big<0\big>^{j-r})$$
if $\Lambda,\Omega$ are even integral, and set $\widetilde A_j(\Lambda,\Omega)=0$
otherwise.
Then the $\Lambda$th coefficient of
$f|\widetilde T_j(p^2)$ is
$$\sum_{p\Lambda\subseteq\Omega\subseteq{1\over p}\Lambda}
\widetilde A_j(\Lambda,\Omega) c(\Omega).$$
\endproclaim

\demo{Proof}  By Theorem 2.4, the $\Lambda$th coefficient of $f|\widetilde T_j(p^2)$
is $$\sum_{p\Lambda\subseteq\Omega\subseteq\frac{1}{p}\Lambda}
\widetilde A_j(\Lambda,\Omega)c(\Omega)$$
where
$$\widetilde A_j(\Lambda,\Omega)=\chi'(p^{j-r})p^{E'_j(\Lambda,\Omega)}
\sum_{\ell=0}^j\sum_{{\cls U}\atop{\dim U=\ell-r}}
\beta(n-\ell,j-\ell) R^*(V,U)\widetilde\G(U)$$
with $E'_j(\Lambda,\Omega)=j(k-n)+k(n_2-n_0)+n_0(n-n_2)$
 and $V=\Lambda_1/p\Lambda_1$.
The number of ways to extend a dimension $\ell-r$ subspace $U$ of $V$ to a dimension $j-r$
subspace $W$ of $V$ is 
$$\beta((n-r)-(\ell-r),(j-r)-(\ell-r))=\beta(n-\ell,j-\ell);$$
thus
$$\widetilde A_j(\Lambda,\Omega)
=\chi'(p^{j-r})p^{E'_j(\Lambda,\Omega)}
\sum_{{\cls W}\atop{\dim W=j-r}}R^*(V,W)
\sum_{{\cls U}\atop{\dim U\le j-r}} R^*(W,U)\widetilde\G(U).$$
By Lemma 3.2, we get
$$\align
\widetilde A_j(\Lambda,\Omega)
&=\chi'(p^{j-r})p^{E_j(\Lambda,\Omega)}
\sum_{{\cls W}\atop{\dim W=j-r}}R^*(V,W)
R^*(W\perp\big<2\big>,\big<0\big>^{j-r})\\
&=\chi'(p^{j-r})p^{E_j(\Lambda,\Omega)}
\left( R^*(V,\big<0\big>^{j-r})+2 R^*(V,\big<0\big>^{j-r-1}\perp\big<-2\big>)\right)\\
&=\chi'(p^{j-r})p^{E_j(\Lambda,\Omega)}
R^*(V\perp\big<2\big>,\big<0\big>^{j-r}),
\endalign$$
proving the theorem.  $\square$
\enddemo

\bigskip
\head{\bf \S4. Hecke operators on Siegel theta series of weight $k+1/2$}\endhead
\smallskip

Throughout this section, we assume $L$ is a rank $2k+1$ lattice with an even
integral, positive definite quadratic form $Q$; we fix $n\le 2k+1$.
As in the case when rank $L$ is even
([13], [14]), we prove a generalised Eichler Commutation
Relation, and from this show the average theta series is an eigenform for 
$\widetilde T_j(p^2)$ ($p\nmid$ level $L$), computing the eigenvalues.
As in Theorem 1.1,
we use $B_Q$ to denote the symmetric bilinear form associated to $Q$,
$N$ the level of
$L$, $\chi$ the character of $\theta^{(n)}(L)$, and $\chi'$ the character
defined by $\chi'(d)=\chi(d)\left({(-1)^{k+1}\over |d|}\right) (\sgn\,d)^{k+1}.$
With $L=\Z v_1\oplus\cdots\oplus\Z v_{2k+1}$, the matrix for $Q$ relative to this
basis for $L$ is given by $Q=\left(B_Q(v_h,v_i)\right).$  Thus for $C\in\Z^{2k+1,n}$
and with $(x_1\,\ldots\,x_n)=(v_1\,\ldots\,v_{2k+1})C$,
$^tCQC=\left(B_Q(x_h,x_i)\right)$ 
is the matrix for the quadratic form $Q$
restricted to the (external) direct sum
$\Lambda=\Z x_1\oplus\cdots\oplus\Z x_n.$
Hence $$\theta^{(n)}(L;\tau)=\sum_{x_1,\ldots,x_n\in L} \e\{(B_Q(x_h,x_i))\tau\}.$$

Note that as a sublattice of $L$, we may have $d=\rank(\Z x_1+\ldots+\Z x_n) <n$.
In such a case there exists some $G\in GL_n(\Z)$ so that
$$(x_1\,\ldots\, x_n)G=(x_1'\,\ldots\, x_d'\,0\,\ldots\,0).$$
Still, we can consider $\Lambda$ as a sublattice of $L$ with ``formal rank'' $n$.
Given a sublattice $\Lambda'=\Z x_1'+\cdots+\Z x_d'$ of $L$ with $\rank \Lambda'=d$
and $T'=\big(B_Q(x_h',x_i')\big)$ (a $d\times d$ matrix),
$$\sum_{{x_1,\ldots,x_n\in L}\atop{\Z x_1+\cdots+\Z x_n=\Lambda'}}
\e\left\{\big(B_Q(x_h,x_i)\big)\tau\right\}
=\sum_G\e\left\{\,^tG\pmatrix T'\\&0_{n-d}\endpmatrix G\tau\right\}$$
where $G$ varies over
$$\left\{\pmatrix I_d&0\\*&*\endpmatrix\in GL_n(\Z)\right\}\Big\backslash GL_n(\Z).$$
Thus with $x_{d+1}'=\cdots=x_n'=0$, $\Lambda=\Z x_1'\oplus\cdots\oplus\Z x_n'$
(the external direct sum), we define
$$\e\{\Lambda\tau\}=\sum_G \e\left\{\,^tG\pmatrix T'\\&0_{n-d}\endpmatrix G\tau\right\}$$
where $G$ varies as above.  Then
$$\theta^{(n)}(L;\tau)=\sum_{\Lambda\subseteq L}\e\{\Lambda\tau\}$$
where $\Lambda$ varies over all distinct sublattices of $L$ with formal rank $n$.
(When $x_i,y_i\in L$, we say
$\Z x_1\oplus\cdots\oplus\Z x_n$ and $\Z y_1\oplus\cdots\oplus\Z y_n$ are 
distinct sublattices
of $L$ with formal rank $n$ when
$\Z x_1+\cdots+\Z x_n\not=\Z y_1+\cdots+\Z y_n.$)

\smallskip\noindent
{\bf Remark.} For $x_i\in L$, $\Lambda=\Z x_1\oplus\cdots\oplus\Z x_n$, and 
$\Lambda'=\Z x_1+\cdots+\Z x_n$, 
we have $\e\{\Lambda\tau\}=o(\Lambda')\e^*\{\Lambda\tau\}$
since, with $d=\rank\Lambda'$ (as a sublattice of $L$),
$$O(\Lambda'\perp\big<0\big>^{n-d})=
\left\{\pmatrix E'&0\\*&*\endpmatrix\in GL_n(\Z):\ E'\in O(\Lambda')\ \right\}.$$

\proclaim{Proposition 4.1}  For $p$ a prime not dividing $N$ and $1\le
j\le n$, 
take $\Omega\subseteq\frac{1}{p}L$ so that $\Omega$ is even integral and has
formal rank $n$;
decompose $\Omega$ as ${1\over p}\Omega_0\oplus\Omega_1\oplus p\Omega_2$
where $\Omega_i\subseteq L$ and $\Omega_0\oplus\Omega_1$ is primitive in $L$ modulo $p$, meaning
that the formal rank of $\Omega_0\oplus\Omega_1$ is its rank in $L$, which is also
the dimension of $\overline{\Omega_0\oplus\Omega_1}$ in $L/pL$.  
Let $r_i$ be the (formal) rank
of $\Omega_i$; set
$$E(\ell,t,\Omega)=t(k-n)+t(t-1)/2+\ell(k-r_0-r_1)+\ell(\ell-1)/2,$$
and set
$$\align
\widetilde c_j(\Omega)
&= \sum_{\ell,t}p^{E(\ell,t,\Omega)} R^*(\Omega_1/p\Omega_1\perp\big<2\big>,\big<0\big>^{\ell})\\
&\qquad \cdot
\delta(k-r_0-\ell,t)\beta(r_2,t)\beta(n-r_0-\ell-t,n-j);
\endalign$$
if $\chi'(p)=1$, and
$$\align
\widetilde c_j(\Omega)
&= \sum_{\ell,t}(-1)^{\ell}
p^{E(\ell,t,\Omega)} R^*(\Omega_1/p\Omega_1\perp\big<2\big>,\big<0\big>^{\ell})\\
&\qquad \cdot
\mu(k-r_0-\ell,t)\beta(r_2,t)\beta(n-r_0-\ell-t,n-j)
\endalign$$
if $\chi'(p)=-1$.
Then
$$\theta^{(n)}(L;\tau)|\widetilde T_j(p^2)
=\sum_{\Omega}\widetilde c_j(\Omega)\e\{\Omega\tau\}$$
where $\Omega$ varies over all even integral sublattices of ${1\over p}L$ that have
(formal) rank $n$. 
\endproclaim

\demo{Proof} 
The proof is virtually identical to that of 
Proposition 1.4 of [13] (see also Proposition 2.1 of [14]),
so here we merely give an indication of how this is done.

By the definitions of $T_j(p^2)$ and $\widetilde T_j(p^2)$, we have
$$\theta^{(n)}(L;\tau)|\widetilde T_j(p^2)
=\sum_{{\Lambda\subseteq L}\atop{p\Lambda\subseteq\Omega\subseteq{1\over p}\Lambda}}
\widetilde A_j(\Omega,\Lambda)\e\{\Omega\tau\}$$
where $\widetilde A_j(\Omega,\Lambda)$ is defined in Theorem 3.3.
(Note that at the end of the proof of Theorem 2.4 we made a change of variables that
we do not make here.)  Since $p\not=2$ and $\Omega\subseteq{1\over
p}L$, $\Omega$ is even integral exactly when it is integral, so
$\widetilde A_j(\Omega,\Lambda)=0$ when $\Omega$ is not integral.
 Interchanging the order of summation, we have
$$\theta^{(n)}(L;\tau)|\widetilde T_j(p^2)
=\sum_{{\Omega\subseteq {1\over p}L}\atop{\Omega\text{ integral}}}
\sum_{p\Omega\subseteq\Lambda\subseteq({1\over p}\Omega\cap L)}
\widetilde A_j(\Omega,\Lambda)\e\{\Omega\tau\}.$$
So one fixes $\Omega$, and follows the procedure of [13] to construct all the
$\Lambda$ in the inner sum, keeping track of the data carried in
$\widetilde A_j(\Omega,\Lambda)$. 
$\square$
\enddemo

\proclaim{Proposition 4.2}  Let $L$ be as above, and fix a prime $p\nmid N$; choose
$j$ so that $1\le j\le n$ and $j\le k$.  
We say a lattice $K$ is a $p^j$-neighbour of $L$ if $K\in\gen L$ and
$$L=L_0\oplus L_1\oplus L_2,\quad K={1\over p}L_0\oplus L_1\oplus pL_2$$
with $\rank L_0=\rank L_2=j$.
\itemitem{(a)} The number of $p^j$-neighbours of $L$ is $p^{j(j-1)/2}\beta\delta(k,j).$
\itemitem{(b)} 
Take $\Omega\subseteq\frac{1}{p}L$ so that $\Omega$ is even integral and has formal
rank $n$; decompose $\Omega$ as in Proposition 4.1.  Set 
$$\align
b_j(\Omega)&=p^{(j-r_0)(j-r_0-1)/2}\sum_{\ell=0}^{j-r_0}
p^{\ell(k-j-r_1+\ell)}
R^*(\Omega_1/p\Omega_1\perp\big<2\big>,\big<0\big>^{\ell})\\
&\qquad \cdot
\delta(k-r_0-\ell,j-r_0-\ell)\beta(k-r_0-r_1,j-r_0-\ell)
\endalign$$
if $\chi'(p)=1$,
$$\align
b_j(\Omega)&=p^{(j-r_0)(j-r_0-1)/2}\sum_{\ell=0}^{j-r_0}
(-1)^{\ell} p^{\ell(k-j-r_1+\ell)}
R^*(\Omega_1/p\Omega_1\perp\big<2\big>,\big<0\big>^{\ell})\\
&\qquad \cdot
\beta(k-r_0-\ell,j-r_0-\ell)\delta(k-r_0-r_1,j-r_0-\ell)
\endalign$$
if $\chi'(p)=-1$.  Then
$$\sum_{K_j}\theta^{(n)}(K_j;\tau)=\sum_{\Omega}b_j(\Omega)\e\{\Omega\}$$
where $K_j$ varies over all $p^j$-neighbours of $L$, and $\Omega$ varies over all even
integral sublattices of ${1\over p}L$ with (formal) rank $n$. 
\endproclaim

\demo{Proof}  
The proof is virtually identical to that of Proposition 1.5 of [13] (see also 
Proposition 2.2 of [14]), so here
we merely give an indication of the proof.

To construct all $p^j$-neighbours $K$ of $L$, we begin by choosing
a $j$-dimensional totally isotropic subspace $\overline C$ of $L/pL$.
Thus $L=C\oplus D\oplus J$ where $\overline C\oplus\overline D\simeq\H^j$,
$B_Q(C\oplus D,J)\equiv0\ (\mod p).$  Set $K'=C\oplus pD\oplus pJ$
(the preimage of $\overline C$); then
in $K'/pK'$ (scaled by $1/p$), $\overline C\oplus \overline{pD}\simeq\H^j$
with $\overline{pD}$ totally isotropic.  Thus we can refine $\overline C$ to
a totally isotropic subspace $\overline{C'}$; set 
$K=\frac{1}{p}C'\oplus pD\oplus J$ (the preimage of $(\overline{C'})^{\perp}$).
Given $\Omega\subseteq\frac{1}{p}L$ with
$\Omega$ decomposed as described above, $\Omega$ lies
in $K$ exactly when $\overline\Omega_0\oplus\overline\Omega_1\subseteq\overline C^{\perp}$ in $L/pL$,
and $\overline\Omega_0\subseteq\overline{C'}$ in $K'/pK'$, and we can easily
count how often this is the case.  $\square$
\enddemo

Next we have a generalised Eichler Commutation Relation.

\proclaim{Theorem 4.3}  Let
$L$ be a lattice of rank $2k+1$ equipped with a positive definite
quadratic form $Q$ of level $N$, and $1\le n\le 2k+1$;
fix a prime $p$ so that $p\nmid N$.
Take $j$ so that $1\le j\le n$ and $j\le k$; for $0\le q\le j$, set
$$\align
u_q(j)&=(-1)^qp^{q(q-1)/2}\beta(n-j+q,q),
\qquad T_j'(p^2)=\sum_{q=0}^ju_q(j)\widetilde T_{j-q}(p^2),\\
v_q(j)&=\cases
(-1)^q\beta(k-n+q-1,q)\delta(k-j+q,q)&\text{if $\chi'(p)=1$,}\\
(-1)^q\delta(k-n+q-1,q)\beta(k-j+q,q)&\text{if $\chi'(p)=-1$.}
\endcases
\endalign$$
Then
$$\theta^{(n)}(L)|T'_j(p^2)=\sum_{q=0}^j v_q(j)\left(
\sum_{K_{j-q}} \theta^{(n)}(K_{j-q})\right)$$
where $K_{j-q}$ runs over all $p^{j-q}$-neighbours of $L$ (as defined in
Proposition 6.3).
\endproclaim

\demo{Proof}  
Fix $\Omega\subseteq\frac{1}{p}L$ as in Propositions 4.1 and 4.2;
suppose that $\chi'(p)=1$.  Then with $\overline\Omega_1=\Omega_1/p\Omega_1$,
$$\align
\sum_{q=0}^j u_q(j)\widetilde c_{j-q}(\Omega)
&=\sum_{\ell,t} p^{E(\ell,t,\Omega)} 
R^*(\overline{\Omega_1}\perp\big<2\big>,\big<0\big>^{\ell})
\delta(k-r_0-\ell,t)\beta(r_2,t)\\
&\qquad \cdot
\sum_{q=0}^j u_q(j) \beta(n-r_0-\ell-t,n-j+q)
\endalign$$
where $E(\ell,t,\Omega)$ is defined in Proposition 4.1.
Using the identities
$$\align
\beta(m+r,r+q)\beta(r+q,q)
&={\mu(m+r,r)\mu(m,q)\mu(r+q,q)\over \mu(r+q,q)\mu(r,r)\mu(q,q)}\\
&=\beta(m+r,r)\beta(m,q)
\endalign$$
and $\beta(m,q)=p^q\beta(m-1,q)+\beta(m-1,q-1)$, we find that
when $m=j-r_0-\ell-t\ge 1$, 
$\sum_{q=0}^m u_q(j) \beta(m,n-j+q)=0.$
Thus
$$\align
\sum_q u_q(j)\widetilde c_{j-q}(\Omega)
&=\sum_{\ell}p^{E(\ell,j-r_0-\ell,\Omega)} 
R^*(\overline\Omega_1\perp\big<2\big>,\big<0\big>^{\ell})\\
&\qquad \cdot
\delta(k-r_0-\ell,j-r_0-\ell) \beta(r_2,j-r_0-\ell).
\endalign$$

On the other hand,
$$\align
\sum_q v_q(j) b_{j-q}(\Omega)
&=p^{(j-r_0)(j-r_0-1)/2}
\sum_{\ell} p^{\ell(k-r_1+\ell-j)} 
R^*(\overline\Omega_1\perp\big<2\big>,\big<0\big>^{\ell}\\
&\qquad \cdot
{\delta(k-r_0-\ell,j-r_0-\ell)\over \mu(j-r_0-\ell,j-r_0-\ell)}
S(j-r_0-\ell)
\endalign$$
where
$$S(m)=\sum_{q=0}^m (-1)^q p^{q(q+1)/2-qm}\mu(k-r_0-r_1,m-q)
\mu(k-n-1+q,q)\beta(m,q).$$
Using that $n=r_0+r_1+r_2$ and that
$\beta(m,q)=\beta(m-1,q)+p^{m-q}\beta(m-1,q-1),$ we find
that 
$$S(m)=p^{k-n}(p^{r_2-m+1}-1)S(m-1)=p^{m(k-n)}\mu(r_2,m).$$
Thus $$\sum_qv_q(j) b_{j-q}(\Omega)=\sum_qu_q(j)\widetilde c_{j-q}(\Omega).$$

The case when $\chi'(p)=-1$ is virtually the same, and so the details are
left to the reader.  $\square$
\enddemo

In the next corollary, we average across the generalised Eichler Commutation
Relation to show that $\theta^{(n)}(\gen L)$ is a Hecke eigenform for primes
$p\nmid N$, where
$$\theta^{(n)}(\gen L)=\sum_{\cls K\in\gen L}{1\over o(K)}\theta^{(n)}(K)$$
(here $\cls K$ varies over all isometry classes within the genus of $L$).

\proclaim{Corollary 4.4} With $p$ a prime, $p\nmid N$, and $1\le j\le n$
with $j\le k$,
$$\theta^{(n)}(\gen L)|T'_j(p^2)=\lambda_j(p^2) \theta^{(n)}(\gen L)$$
where
$$\lambda_j(p^2)=\cases
p^{j(j-1)/2+j(k-n)}\beta(n,j)\delta(k,j)&\text{if $\chi'(p)=1$,}\\
p^{j(j-1)/2+j(k-n)}\beta(n,j)\mu(k,j)&\text{if $\chi'(p)=-1$.}
\endcases$$
\endproclaim

\demo{Proof} First note that for $K\in\gen L$, we have $\disc K=\disc L$, so
$K$ is a $p^m$-neighbour of $L$ (as defined in Proposition 4.2)
if and only if $pL\subseteq K\subseteq{1\over p}L$, and either $\mult_{\{L:K\}}(p)=m$
or $\mult_{\{L:K\}}(1/p)=m$ where $\{L:K\}$ denotes the invariant
factors of $K$ in $L$.  
Classifying the $p^m$-neighbours into isometry classes,
we see that the number of $p^m$-neighbours of $L$ in $\cls K\in\gen L$ is
$${\#\{{\text{isometries }}\sigma:\ pL\subseteq\sigma K\subseteq{1\over p}L,
\ \mult_{\{L:\sigma K\}}(p)=m\ \}\over o(K)}$$
(since $\sigma K=\sigma' K$ if and only if $\sigma^{-1}\sigma'\in O(K)$). Also,
using Proposition 4.2 (a), 
$$\align
&\sum_{\cls L'\in\gen L}
{\#\{{\text{isometries }}\sigma:\ pL'\subseteq\sigma K\subseteq{1\over p}L',
\ \mult_{\{L':\sigma K\}}(p)=m\ \}\over o(L') o(K)}\\
&\qquad = {1\over o(K)} \sum_{\cls L'}
{\#\{{\text{isometries }}\sigma:\ pK\subseteq\sigma L'\subseteq{1\over p}K,
\ \mult_{\{K:\sigma L'\}}(p)=m\ \}\over o(L')}\\
&\qquad = {1\over o(K)}\#\{ p^m{\text{-neighbours of }}K\}\\
&\qquad = {1\over o(K)} p^{m(m-1)/2}\beta\delta(k,m).
\endalign$$
Thus
$$\theta^{(n)}(\gen L)|T'_j(p^2)=\lambda_j(p^2) \theta^{(n)}(\gen L)$$
where
$$\lambda_j(p^2)=\sum_{q=0}^j v_q(j) p^{(j-q)(j-q-1)/2}\beta\delta(k,j-q).$$
When $\chi'(p)=1$, 
$\lambda_j(p^2)=p^{j(j-1)/2}\frac{\delta(k,j)}{\mu(j,j)}S(j)$ where
$$S(j)=\sum_{q=0}^j (-1)^q p^{q(q+1)/2-qj} \beta(j,q)\mu (k-n+q-1,q)\mu(k,j-q).$$
Using the identity $\beta(j,q)=\beta(j-1,q)+p^{j-q}\beta(j-1,q-1),$ we find that for
$1\le d\le j$, 
$$S(j)=p^{d(k-n)}\mu(n-j+d,d) S(j-d) = p^{j(k-n)} \mu(n,j),$$
proving the corollary when $\chi'(p)=1$.

The case $\chi'(p)=-1$ is virtually identical. $\square$
\enddemo

Just as in the integral weight case (\S3 of [14]), we have the following.

\proclaim{Theorem 4.5} When $1\le a\le n-k$,
$p$ is a prime not dividing $N$, and $T'_j(p^2)$ is defined as in Theorem 4.3,
$$\theta^{(n)}(L)|T'_{k+a}(p^2)=0.$$
\endproclaim

\demo{Proof}  One argues exactly as in the integral weight case (see \S 3 [14]):
First, using Proposition 4.1, one shows
$$\widetilde c_{k+a}(\Omega)=\sum_{q=0}^k w_q(a)\widetilde c_{k-q}(\Omega)$$ where
$w_q(a)=(-1)^q p^{q(q+1)/2}\beta(a+q-1,q)\beta(n-k+q,a+q).$
Hence 
$$\theta^{(n)}(L)|\widetilde T_{k+a}(p^2)
=\theta^{(n)}(L)|\sum_{q=0}^k w_q(a)\widetilde T_{k-q}(p^2).$$
Then one shows 
$$\sum_{q=0}^r\beta(n-q,r-q) T'_q(p^2)=\widetilde T_r(p^2).$$
Substituting and using induction on $a$ yields the result. $\square$
\enddemo

\bigskip
\head{\S5. Bounding Hecke eigenvalues}\endhead
\smallskip

In [16] we used the formula from [6] to bound the eigenvalues of Hecke operators
acting on integral weight Siegel forms.  Using Theorem 3.3, we can do
this for half-integral weight.  We obtain the following bound.

\proclaim{Theorem 5.1}  Let $f$ be a nonsingular Siegel modular form of degree $n$
and weight $k+1/2$ so that for some $j$, $1\le j\le n$,
$f|\widetilde T_j(p^2)=\lambda_j(p^2)\,f$ for almost all primes $p$.
Also suppose that the Fourier coefficients $c(\Lambda)$ of $f$ satisfy the bound
$|c(\Lambda)|\ll_f(\disc\Lambda)^{k/2+1/4-\gamma}$.
(Note that with $\gamma=0$, this is the trivial bound for cusp forms.)
Then
$$|\lambda_j(p^2)|\ll_f\,p^M \text{ (as $p\mapsto\infty$)}$$
where $M=\frac{1}{4}(j+n-2\gamma+\frac{1}{2})^2+\frac{1}{6}(j-n+2\gamma-1)^2+j(k-n).$
When $\gamma=0$, we can take $M=\frac{1}{4}(j+n+\frac{1}{2})^2+j(k-n).$
\endproclaim

\demo{Proof}  The proof is somewhat similar to that of Proposition 4.1, and almost
identical to that in [16]:  Given $\Lambda$, we construct all the even integral
$\Omega$ so that $p\Lambda\subseteq\Omega\subseteq\frac{1}{p}\Lambda$, keeping
track of the relevant data for the formula of Theorem 3.3.  Note that since we are
bounding the $\Lambda$th coefficient of $f|\widetilde T_j(p^2)$, we can replace
$\chi'(p^{j-r})$ by 1 in this formula.

Since $f$ is nonzero, we can choose $\Lambda$ so that $c(\Lambda)\not=0$ and
$\disc\Lambda\not=0$.  Also, assume $p\nmid 2\disc\Lambda$.  We first partition the lattices
$\Omega$, $p\Lambda\subseteq\Omega\subseteq\frac{1}{p}\Lambda$, according to the
invariant factors $\{\Lambda:\Omega\}$.  Choosing $n_0,n_2\ge 0$ so that
$n_0+n_2\le j$, we construct all integral $\Omega$ so that
$$\Omega=p\Lambda_0\oplus\Lambda_1\oplus\frac{1}{p}\Lambda_2
\text{ where }\Lambda=\Lambda_0\oplus\Lambda_1\oplus\Lambda_2,$$
$n_i=\rank\Lambda_i$.  In this process, we compute
$R^*(\Lambda_1/p\Lambda_1\big<2\big>,\big<0\big>^{j-r})$
where $r=n_0+n_2$.
(Note that since $p\nmid\disc\Lambda$, $\Lambda/p\Lambda$ is a regular quadratic space;
so when $\Omega$ is integral, we must have that, in $\Lambda/p\Lambda$,
$\overline\Lambda_2$ is totally isotropic and $\overline\Lambda_1$ is orthogonal
to $\overline\Lambda_2$.  Consequently we must have $n_2\le n_0$.)

First we choose $\overline\Delta_2$ to be a totally isotropic, dimension $n_2$ subspace
of $\Lambda/p\Lambda$.  Since $\Lambda/p\Lambda$ is regular, the lemma of \S4 of [16]
tells us the number of choices for $\overline\Delta_2$ is bounded by
$$4^{n_2}p^{n_2(n-n_2)-n_2(n_2+1)/2}.$$
Note that $\Lambda/p\Lambda=(\overline\Delta_2\oplus\overline\Delta'_2)\perp J$
where $\overline\Delta_2\oplus\overline\Delta'_2\simeq\H^{n_2}$ and $J$ is regular;
thus $\overline\Delta^{\perp}_2=\overline\Delta_2\perp J$.

Next we extend $\overline\Delta_2$ to $\overline\Delta_2\oplus W$ where
$W\subseteq\overline\Delta^{\perp}_2$ and either $W\simeq\big<0\big>^{j-r}$
or $W\simeq\big<0\big>^{j-r-1}\perp\big<-2\big>$; when 
$W\simeq\big<0\big>^{j-r-1}\perp\big<-2\big>$, we count $W$ with multiplicity 2.
So the number of such $W$, counted with the appropriate multiplicity, is
$R^*(J\perp\big<2\big>,\big<0\big>^{j-r};$ by \S 4 of [16], this is bounded by
$$4^{j-r}p^{(j-r)(n-2n_2+1-j+r)-(j-r)(j-r+1)/2}.$$

Now we extend $\overline\Delta_2\oplus W$ to $\overline\Delta_2\oplus\overline\Lambda_1$
where $\overline\Lambda_1\subseteq\overline\Delta^{\perp}_2$.  The number of choices
for $\overline\Lambda_1$ is $\beta(n-j+n_0-n_2,n-j)$, which by \S4 of [16] is bounded by
$$2^{n-j}p^{(n-j)(n_0-n_2)}.$$

Now let $\Delta$ be the preimage in $\Lambda$ of $\overline\Delta_2$, and $\Omega'$
the preimage in $\Lambda$ of $\overline\Delta_2\oplus\overline\Lambda_1$.  Thus
$$\align
\Lambda&=\Lambda_0\oplus\Lambda_1\oplus\Delta_2,\\
\Delta&=p\Lambda_0\oplus p\Lambda_1\oplus\Delta_2,\\
\Omega'&=p\Lambda_0\oplus\Lambda_1\oplus\Delta_2.
\endalign$$
So $Q(\Delta_2)\equiv0\ (\mod p)$, $B_Q(\Delta_2,\Lambda_1)\equiv0\ (\mod p)$
where $Q$ denotes the quadratic form on $\Lambda$, $B_Q$ the corresponding
symmetric bilinear form.

Our final step in constructing $\Omega$ is to refine our choice of $\Delta_2$
so that $Q(\Delta_2)\equiv 0\ (\mod p^2)$.
For this, we work in $\Delta/p\Delta$, with the quadratic form scaled by $1/p$.
We extend $\overline{p\Omega'}=\overline{p\Lambda}_1$ to
$\overline{p\Lambda}_1\oplus\overline\Lambda_2$ where $\overline\Lambda_2$ is 
totally isotropic of dimension $n_2$ and independent of $\overline{p\Lambda}$.
As discussed in [16], the number of choices for $\overline\Lambda_2$ is bounded
by $p^{n_2(n_0-n_2)}.$  Now we take $p\Omega$ to be the preimage in $\Delta$
of $\overline{p\Lambda}_1\oplus\overline\Lambda_2$.
Note that $\disc\Omega=p^{2(n_0-n_2)}\disc\Lambda$.

Thus, using Theorem 3.3 and the assumed bound on the Fourier coefficients of $f$,
we see the $\Lambda$th coefficient of $f|\widetilde T_j(p^2)$ is bounded by
$$\sum_{n_0+n_2\le j} 2^{n+j-2n_0}p^{E(n_0,n_2)}(\disc\Lambda)^{k/2+1/4-\gamma}$$
where $E(n_0,n_2)=-n_0^2+n_0(j+n-2\gamma+1/2)-\frac{3}{2}n_2^2+n_2(j-n+2\gamma-1)
+j(k-n).$  $E(n_0,n_2)$ has its maximum when $n_0=\frac{1}{2}(j+n-2\gamma+1/2)$,
$n_2=\frac{1}{3}(j-n+2\gamma-1)$.  Since $n_0,n_2$ are actually restricted to
be non-negative, when $\gamma=0$ we see that $E(n_0,n_2)$ has its maximum when
$n_0=\frac{1}{2}(j+n+1/2)$, $n_2=0$.  With $M$ the maximum value of $E(n_0,n_2)$,
we have
$$|\lambda_j(p^2)c(\Lambda)|\ll_f\,(\disc\Lambda)^{k/2+1/4-\gamma}p^M,$$
and hence $|\lambda_j(p^2)|\ll_{f,M}\,p^M.$ $\square$
\enddemo

\bigskip
\head{\S6. Hecke-stability of the Kitaoka space}\endhead
\smallskip

In [8], Kitaoka identified a subspace of integral weight Siegel modular forms
$f$ where $c_f(\Lambda)=c_f(\Lambda')$ whenever $\Lambda'\in\gen\Lambda$, and he
showed the space is stable under the Hecke operators associated to primes not
dividing the level.  Here we show that the half-integral weight Kitaoka subspace
is invariant under all Hecke operators; note that our argument also works for the 
integral weight Kitaoka subspace (without the restriction that the
prime does not divide the level).

\proclaim{Proposition 6.1}  For $p$ a prime
and $\Lambda,\Lambda'$ lattices in the same
genus, there is a bijective map
$$\sigma:\left\{\Omega:\ p\Lambda\subseteq\Omega\subseteq{1\over p}\Lambda\ \right\}
\to \left\{\Omega':\ p\Lambda'\subseteq\Omega'\subseteq{1\over p}\Lambda'\ \right\}$$
so that $\sigma(\Omega)\in\gen\Omega$.  Further, with
$p\Lambda\subseteq\Omega\subseteq{1\over p}\Lambda$, we have
$$\{\Lambda:\Omega\}=\{\Lambda':\sigma(\Omega)\}
\text{ and }
\alpha_j(\Lambda,\Omega)=\alpha_j(\Lambda',\sigma(\Omega))$$
where $\{\Lambda:\Omega\}$ denotes the invariant factors of $\Omega$
in $\Lambda$.
Hence $\sigma$ also gives a bijection between 
$$\left\{\Omega:\ p\Lambda\subseteq\Omega\subseteq\Lambda\ \right\} \text{ and }
\left\{\Omega':\ p\Lambda'\subseteq\Omega'\subseteq\Lambda'\ \right\}.$$
\endproclaim

\demo{Proof}  Since $\Lambda,\Lambda'$ are in the same genus, we can choose
an isometry $\sigma_p:\Q_p\Lambda\to\Q_p\Lambda'$ so that
$\sigma_p(\Z_p\Lambda)=\Z_p\Lambda'$.  Then for
$\Omega\subseteq{1\over p}\Lambda$,
and working in $\Q_p\Lambda'$, we set
$$\sigma(\Omega)=\sigma_p(\Z_p\Omega)\cap{1\over p}\Lambda'.$$
(So $\sigma(\Lambda)=\Lambda'$.)

Now take $\Omega$ so that $p\Lambda\subseteq\Omega\subseteq{1\over p}\Lambda$;
thus
$$p\Lambda'=\sigma(p\Lambda)\subseteq\sigma(\Omega)\subseteq
\sigma\left({1\over p}\Lambda\right)={1\over p}\Lambda'.$$
Also, $$\sigma_p(\Z_p\Omega)\subseteq\sigma_p\left({1\over
p}\Z_p\Lambda\right)
={1\over p}\Z_p\Lambda',$$
so
$$\Z_p\sigma(\Omega)=\sigma_p(\Z_p\Omega)\cap{1\over p}\Z_p\Lambda'
=\sigma_p(\Z_p\Omega)\simeq\Z_p\Omega.
$$
For all primes $q\not=p$, since
$p\Lambda'\subseteq\sigma(\Omega)\subseteq{1\over p}\Lambda',$
we have
$$\Z_q\Lambda'=p\Z_q\Lambda'\subseteq\Z_q\sigma(\Omega)\subseteq{1\over
p}\Z_q\Lambda'
=\Z_q\Lambda';$$
thus $\Z_q\sigma(\Omega)=\Z_q\Lambda'$.  Similarly,
$\Z_q\Omega=\Z_q\Lambda$,
so 
$$\Z_q\sigma(\Omega)=\Z_q\Lambda'\simeq\Z_q\Lambda=\Z_q\Omega.$$
Hence $\Omega'\in\gen\Omega.$

To see $\sigma$ is a bijection as claimed, take
$\Omega'\subseteq{1\over p}\Lambda'$;
working in $\Q_p\Lambda$, set
$$\sigma'(\Omega')=\sigma_p^{-1}(\Z_p\Omega')\cap{1\over p}\Lambda.$$
Then one easily checks that $\sigma\circ\sigma'(\Omega')=\Omega'$, and
for $\Omega\subseteq\frac{1}{p}\Lambda$, $\sigma'\circ\sigma(\Omega)=\Omega$.

Next, we show that for $p\Lambda\subseteq\Omega\subseteq{1\over p}\Lambda$,
we have
$\{\Lambda:\Omega\}=\{\Lambda':\sigma(\Omega)\}$.
First recall that for all
primes $q\not=p$, $\Z_q\Omega=\Z_q\Lambda$ and 
$\Z_q\sigma(\Omega)=\Z_q\Lambda'$.  Hence
$$\align
\{\Lambda':\sigma(\Omega)\}&=\{\Z_p\Lambda':\Z_p\sigma(\Omega)\}
=\{\sigma_p(\Z_p\Lambda):\sigma_p(\Z_p\Omega)\}\\
&=\{\Z_p\Lambda:\Z_p\Omega\}
=\{\Lambda:\Omega\}.
\endalign$$

Finally, we show that with $p\Lambda\subseteq\Omega\subseteq{1\over p}\Lambda$,
$\alpha_j(\Lambda,\Omega)=\alpha_j(\Lambda',\sigma(\Omega)).$
Note that since $\Z/p\Z\approx\Z_p/p\Z_p$, the $\Z/p\Z$-space
$(\Lambda\cap\Omega)/p(\Lambda+\Omega)$ and the
$\Z_p/p\Z_p$-space $(\Z_p\Lambda\cap\Z_p\Omega)/p(\Z_p\Lambda+\Z_p\Omega)$
can be viewed as isometric (over $\Z_p/p\Z_p$).  Thus
$$\align
\alpha_j(\Lambda,\Omega)&=\alpha_j(\Z_p\Lambda,\Z_p\Omega)\\
&=\alpha_j\big(\sigma_p(\Z_p\Lambda),\sigma_p(\Z_p\Omega)\big)
=\alpha_j(\Z_p\Lambda',\Z_p\sigma(\Omega))
=\alpha_j(\Lambda',\sigma(\Omega)).
\endalign$$
This proves the proposition. $\square$
\enddemo

Note that when $\Lambda,\Lambda'$ are oriented, $\Omega$ inherits its orientation
from $\Lambda$, and $\sigma(\Omega)$ inherits its orientation from $\sigma(\Lambda)=\Lambda'$.
Then the aboves proposition together with Theorem 3.4  immediately
gives us the following.

\proclaim{Theorem 6.2}  
For $k,n,N\in\Z_+$ and $\chi$ a character modulo $N$, the subspace
$$\left\{f\in\M_{k+1/2}(\Gamma_0^{(n)}(N),\chi)):\ c_f(\Lambda)=c_f(\Lambda')
\text{ when }\Lambda'\in\gen\Lambda\ \right\}$$
is stable under the full Hecke-algebra.
\endproclaim

\bigskip
\head{\S7.  A transparent Hecke-correspondence}\endhead
\smallskip

In [15], we developed a formula for the action of Hecke operators on
Jacobi modular forms; here we observe that this is almost identical to our
formula in Theorem 3.4.  From this we easily obtain a Hecke-correspondence.
We begin by introducing notation and terminology for Jacobi modular forms
that we will use in our correspondence.

A Jacobi modular form of index $1$, weight $k+1$, Siegel degree $n+1$, level $N$,
and character $\chi'$ is an analytic function $F:\H_{(n)}\times\C^{1,n}\to \C$
so that
$$\widehat F\pmatrix \tau&^tZ\\Z&\tau'\endpmatrix = F(\tau,Z)\e\{2\tau'\}$$
transforms like a degree $n+1$, weight $k+1$, level $N$, character $\chi'$
Siegel modular form under matrices
$$\underline\gamma\in \Gamma_0^{(n,1)}(N)
=\left\{\pmatrix A&0&B&*\\*&1&*&*\\C&0&D&*\\0&0&0&1\endpmatrix:\ 
\pmatrix A&B\\C&D\endpmatrix\in \Gamma_0^{(n)}(N)\ \right\}.$$
(Here $\tau'$ is a formal variable.)
So with $\underline\tau=\pmatrix\tau&^tZ\\Z&\tau'\endpmatrix$,
$\underline\gamma\in\Gamma_0^{(n,1)}(N)$, we have
$(F(\tau,Z)|\underline\gamma)\e\{2\tau'\}=\widehat
F(\underline\tau)|\underline\gamma;$ also,
$$\widehat F(\underline\tau)=\sum_{\underline T}\widehat c(\underline
T)\e\{\underline T\underline\tau\} \text{ and }
F(\tau,Z)=\sum_{T,R} c(T,R) \e\{T\tau+\,^tRZ\}$$
where $T$ varies over $n\times n$ symmetric matrices and $R$ varies
over $1\times n$ matrices so that 
$\underline T=\pmatrix T&^tR\\R&2\endpmatrix$
is even integral and positive semi-definite, and $\widehat
c(\underline T)=c(T,R).$
We say $\widehat F, F$ are even if $2|R$ whenever $\underline T$ is in the support of $\widehat F$.
Let
$J_{k+1,1}(\Gamma_0^{(n,1)}(N),\chi')$ denote the space of Jacobi modular forms
of weight $k+1$, index $1$, Siegel degree $n+1$, level $N$, and character $\chi'$.
Let $J^{\text{even}}_{k+1,1}(\Gamma_0^{(n,1)}(N),\chi')$ denote the subspace of 
even Jacobi modular forms in $J_{k+1,1}(\Gamma_0^{(n,1)}(N),\chi')$.

For 
$$G\in GL_{n+1,1}^J(\Z)=\left\{\underline G=\pmatrix G&0\\R&1\endpmatrix\in GL_{n+1}(\Z)\right\},$$
$c(\,^t\underline G\,\underline T\,\underline G)
=\chi'(\det\underline G) (\det\underline G)^{k+1}
\,c(\underline T).$
Consequently we can define $\widehat F$ as a Fourier series supported on lattices
as follows.  Fix a rank 1 lattice
$\Delta=\Z w\simeq (2)$, and for $\Lambda\oplus\Delta$ a rank
$n+1$ lattice with quadratic form given by the matrix
$\underline T=\pmatrix T&^tR\\R&2\endpmatrix,$ set
$c_{\Delta}(\Lambda\oplus\Delta)=\widehat c(\underline T).$
Let $$\e^*_{\Delta}\{(\Lambda\oplus\Delta)\underline\tau\}=\sum_{\underline G} 
\e\{^t\underline G\,\underline T\,\underline G\,\underline \tau\}$$
where $\underline G$ varies over
 $O(\underline T)\cap GL_{n+1,1}^J(\Z)\backslash GL_{n+1,1}^J(\Z)$
(or, if $\chi'(-1)(-1)^{k+1}=-1$,
over $O(\underline T)\cap SL_{n+1,1}^J(\Z)\backslash SL_{n+1,1}^J(\Z)$
where  $SL_{n+1,1}^J(\Z)=SL_{n+1}(\Z)\cap GL_{n+1,1}^J(\Z)$).
Then $$\widehat F(\underline\tau)=\sum_{\Lambda\oplus\Delta}
c_{\Delta}(\Lambda\oplus\Delta)\,\e^*_{\Delta}\{(\Lambda\oplus\Delta)\underline\tau\}$$
where $\Lambda\oplus\Delta$ varies over all isometry classes of rank $n+1$, even
integral, positive semi-definite lattices (with $\Delta$ fixed, and 
$\Lambda$ oriented when $\chi'(-1)(-1)^{k+1}=-1$; note that the orientation of $\Delta$
is fixed).
Note that with $G=\pmatrix I&0\\-R&1\endpmatrix$, we have
$\pmatrix G\\&^tG^{-1}\endpmatrix\in\Gamma_0^{(n,1)}(N)$ for any $N$,
and so $c(T,2R)=c(T-2\,^tRR,0)$.
Thus when $F$ is an even Jacobi form, the support of $\widehat F$ is on
lattices $\Lambda'\oplus\Delta\simeq\pmatrix T'&2\,^tR\\2R&2\endpmatrix,$
and we have $\Lambda'\oplus\Delta=\Lambda\perp\Delta$ with
$\Lambda\perp\Delta\simeq\pmatrix T\\&2\endpmatrix,$ $T=T'-2\,^tRR.$

We will see that the space of even Jacobi forms is stable under the
Hecke operators $T_j(p^2)$ provided $p\not=2$; to allow us to handle $p=2$, we
introduce a projection map onto the space of even Jacobi forms.

\proclaim{Proposition 7.1}  For $F\in\J_{k+1,1}(\Gamma_0^{(n,1)}(N),\chi')$
with $4|N$, define
$F|\psi$ by
$$F|\psi=2^{-n}\cdot\sum_{Y\,(2)} F|\pmatrix I_n&&&^tY/2\\&1&Y/2\\&&I_n\\&&&1\endpmatrix$$
where $Y$ varies over $\Z^{1,n}$ modulo 2.  Then $\psi$ maps
$\J_{k+1,1}(\Gamma_0^{(n,1)}(N),\chi')$ onto
$\J_{k+1,1}^{\text{even}}(\Gamma_0^{(n,1)}(N),\chi'),$ and $\psi$ acts as the identity map
on $\J_{k+1,1}^{\text{even}}(\Gamma_0^{(n,1)}(N),\chi')$.
\endproclaim

\demo{Proof} Take $F\in\J_{k+1,1}(\Gamma_0^{(n,1)}(N),\chi')$.
We first show $F|\psi\in\J_{k+1,1}(\Gamma_0^{(n,1)}(N),\chi').$

Take $Y\in\Z^{1,n}$,
$\pmatrix A&B\\C&D\endpmatrix\in \Gamma_0^{(n)}(N)$, and set
$$\delta=\pmatrix I_n&&&^tY/2\\&1&Y/2\\&&I_n\\&&&1\endpmatrix, \quad
\gamma_1=\pmatrix A&&B\\&1&&0\\C&&D\\&0&&1\endpmatrix.$$
Then with $Y'=YD$ and
$$\delta'=\pmatrix I_n&&&^tY'/2\\&1&Y'/2\\&&I_n\\&&&1\endpmatrix,$$
we have $\delta\gamma_1=\gamma_1'\delta'$
where $\gamma_1'\in\Gamma_0^{(n,1)}(N)$ and 
$\chi'(\gamma_1')=\chi'(\gamma_1)$.  Also, since $4|N$, we know
$2\nmid\det D$; thus, as $Y$ varies modulo $2$,
so does $Y'$.  Hence $F|\psi|\gamma_1=\chi'(\gamma_1) F|\psi$.

Next, take $U,V\in\Z^{1,n}$, $w\in\Z$, and set
$$\gamma_2=\pmatrix I&&&^tV\\U&1&V&w\\&&I&-\,^tU\\&&&1\endpmatrix.$$
Then with $\delta$ as above, $\delta\gamma_2=\gamma_2'\delta$ with
$\gamma'_2\in\Gamma_0^{(n,1)}(N)$ and $\chi'(\gamma_2)=\chi'(\gamma_2')$; so
$F|\psi|\gamma_2=\chi'(\gamma_2) F|\psi$.
Since matrices of the form $\gamma_1,\gamma_2$ generate $\Gamma_0^{(n,1)}(N)$,
we have  $F|\psi\in\J_{k+1,1}(\Gamma_0^{(n,1)}(N),\chi').$

Finally, applying $\psi$ to the Fourier expansion for $F(\tau,Z)$, we see
that $F|\psi\in\J_{k+1,1}^{\text{even}}(\Gamma_0^{(n,1)}(N),\chi'),$
and that $\psi$ acts as the identity map on the even Jacobi modular forms.
$\square$
\enddemo

Theorem 3.2 of [15] gives us the following theorem; note that in [15]
we refer to the 
even integral quadratic form on $\Delta$ as the index of the Jacobi form, so index
1 in this paper corresponds to index $(2)$ in [15].

\proclaim{Theorem 7.2}
Take $F\in J_{k+1,1}(\Gamma_0^{(n,1)}(N),\chi')$,
$p$ prime, $1\le j\le n$, and $\Delta\simeq(2)$.
\item{(a)} Suppose $p\nmid N$.  For $U$ a rank $d$ lattice identified
with an even integral quadratic form also denoted by $U$, let
$$\alpha'(U)=\sum_Y \e\{UY/p\}$$
where $Y$ varies over symmetric $d\times d$ matrices modulo $p$ so
that $p\nmid\det Y$; when $d=0$, we take $\alpha'(U)=1$.
With $\Omega_1\oplus\Delta$ an even integral lattice, let
 $R^*_{\Delta}(\overline\Omega_1\oplus\overline\Delta,U)$ denote the number of subspaces
of $\overline\Omega_1\oplus\overline\Delta=\Omega_1/p\Omega_1\oplus \Delta/p\Delta$
isometric to $U$ and independent of $\overline\Delta$.
Given even integral
$\Lambda\oplus\Delta$ and $\Omega\oplus\Delta$ so that
$$\Lambda=\Lambda_0\oplus\Lambda_1\oplus\Lambda_2\oplus\Delta,\qquad
\Omega\oplus\Delta={p}\Lambda_0\oplus\Lambda_1\oplus \frac{1}{p}\Lambda_2\oplus\Delta$$
with $n_i=\rank\Lambda_i$, $r=n_0+n_2$, set
$$\align
A^J_{j,\Delta}(\Lambda,\Omega)&=\chi'(p^{j-r})p^{(k+1)(n_2-n_0)+n_0(n-n_2+2)}\\
&\qquad\cdot
\sum_{{\cls U}\atop{\dim U=j-r}} R^*_{\Delta}(\overline\Omega_1\oplus\overline\Delta,U)
\alpha'(U).
\endalign$$
Then the $\Lambda\oplus\Delta$th
Fourier coefficient of $F|T^J_j(p^2)$ is
$$\sum_{\Omega}A^J_{j,\Delta}(\Lambda,\Omega)\,c_{\Delta}(\Omega\oplus\Delta)$$
where $\Omega$ varies so that
$\Omega\oplus\Delta$ is even integral, and
$p\Lambda'\subseteq\Omega\subseteq{1\over p}\Lambda'$ for some
$\Lambda'$ that satisfies $\Lambda'\oplus\Delta=\Lambda\oplus\Delta$.

\item{(b)} If $p|N$ then the $\Lambda\oplus\Delta$th Fourier coefficient of
$F|T^J_j(p^2)$ is
$$p^{j(n+1-k)}\sum_{\Omega} c_{\Delta}(\Omega\oplus\Delta)$$
where $\Omega\oplus\Delta$ varies so that
$p\Lambda\subseteq\Omega\subseteq\Lambda$
with $[\Lambda:\Omega]=p^j$.
\endproclaim

To complete the incomplete character sum $\alpha'(U)$ when $p\nmid N$, we set
$$\widetilde T^J_j(p^2)=p^{j(k-n-1)}\sum_{0\le\ell\le j}\chi'(p^{j-\ell})
p^{j-\ell}\beta(n-\ell,j-\ell) T^J_{\ell}(p^2).$$
(The coefficients for this
linear combination were not correct in [15], and were later
corrected in an erratum.)

\proclaim{Corollary 7.3} Take $F\in J_{k+1,1}(\Gamma_0^{(n,1)}(N),\chi')$
and $p$ a prime so that $p\nmid N$; using the notation of Theorem 7.2,
set
$$\widetilde A^J_{j,\Delta}(\Lambda,\Omega)
=\chi'(p^{j-r})p^{E_{j,\Delta}(\Lambda,\Omega)}
R^*_{\Delta}(\overline\Omega_1\oplus\overline\Delta,\big<0\big>^{j-r})$$
where $E_{j,\Delta}(\Lambda,\Omega)=j(k-n)+k(n_2-n_0)+n_0(n-n_2)+(j-r)(j-r-1)/2.$
Then the 
$\Lambda\oplus\Delta$th Fourier coefficient of $F|\widetilde T^J_j(p^2)$ is
$$\sum_{\Omega}\widetilde A^J_{j,\Delta}(\Lambda,\Omega) c_{\Delta}(\Omega\oplus\Delta)$$
where $\Omega$ varies as in Theorem 7.2.
\endproclaim

\demo{Proof} To prove this, we need to show that
$$\align
&\sum_{0\le\ell\le j}p^{j-\ell}\beta(n-\ell,j-\ell)
\sum_{{\cls U}\atop{\dim
U=\ell-r}}R^*_{\Delta}(\overline\Omega_1\oplus\overline\Delta,U)
\alpha'(U)\\
&\qquad =
\sum_{r\le\ell\le j}\sum_{{\cls W}\atop{\dim W=j-r}}
R^*_{\Delta}(\overline\Omega_1\oplus\overline\Delta,W)\alpha(W)
\endalign$$
where, identifying $W$ with the even integral matrix for the quadratic form on $W$,
$\alpha(W)$ is the complete character sum
$$\alpha(W)=\sum_{Y\,(p)}\e\{WY/p\},$$
$Y$ varying over all symmetric $(j-r)\times(j-r)$ matrices modulo $p$.
(Note that $\alpha(W)=p^{(j-r)(j-r+1)/2}$ if $W\simeq\big<0\big>^{j-r}$,
and 0 otherwise.)

Given $r\le \ell\le j$ and a dimension $\ell-r$ subspace $U$ of
$\overline\Omega_1\oplus\overline\Delta$ that is independent of $\overline\Delta$,
$p^{j-\ell}\beta(n-\ell,j-\ell)$ is the number of ways to extend $U$ to a
dimension $j-r$ subspace $W$ of
$\overline\Omega_1\oplus\overline\Delta$ that is independent of $\overline\Delta$.
With such $W$, all subspaces of $W$ are necessarily independent of $\overline\Delta$.
Therefore, with $d=j-r$,
$$p^{j-\ell}\beta(n-\ell,j-\ell)
R^*_{\Delta}(\overline\Omega_1\oplus\overline\Delta,U)=
\sum_{{\cls W}\atop{\dim W=d}}
R^*_{\Delta}(\overline\Omega_1\oplus\overline\Delta,W)
R^*(W,U).$$
So we need to show that with $W$ of dimension $d$,
$$\alpha(W)=\sum_{a=0}^d\sum_{{\cls U}\atop{\dim U=a}}
R^*(W,U)\alpha'(U).$$

First suppose $p\not=2$; let $\F=\F_p$.  
Take $\omega\in\F^{\times}$
so that $\left(\frac{\omega}{p}\right)=-1$, and for $a\ge 1$, set
$J_a=\diag\{I_{a-1},\omega\}$.
Note that for $Y\in\F^{d,d}_{\sym}$, $Y$ can be diagonalised and so
either $Y\sim I_a\perp\big<0\big>^{d-a}$ or $Y\sim
J_a\perp\big<0\big>^{d-a}$, some $a$.
In what follows we will
sometimes write $GL_d$ for $GL_d(\F)$.
Then (using notation from \S1),
$$\align
\alpha(W)-1
&=\sum_{a=1}^d\left(
\sum_{Y\sim I_a\perp\big<0\big>^{d-a}}\e\{WY/p\}+
\sum_{Y\sim J_a\perp\big<0\big>^{d-a}}\e\{WY/p\}\right)\\
&=\sum_{a=1}^d \sum_{G\in GL_d} 
{1\over o(I_a\perp\big<0\big>^{d-a})}\e\left\{W\,^tG\pmatrix I_a\\&0\endpmatrix G/p\right\}\\
&\qquad + \sum_{a=1}^d\sum_{G\in GL_d}
{1\over o(J_a\perp\big<0\big>^{d-a})}\e\left\{W\,^tG\pmatrix J_a\\&0\endpmatrix G/p\right\}.
\endalign$$
We have
$$\align
o(I_a\perp\big<0\big>^{d-a})
&=r^*(I_a\perp\big<0\big>^{d-a},I_a)\,r^*(\big<0\big>^{d-a},\big<0\big>^{d-a})\\
&=p^{a(d-a)}\,o(I_a)\,\prod_{i=0}^{d-a-1}(p^{d-a}-p^i),
\endalign$$
and similarly, $o(J_a\perp\big<0\big>^{d-a})
=p^{a(d-a)}\,o(J_a)\,\prod_{i=0}^{d-a-1}(p^{d-a}-p^i).$
Also, writing $GW\,^tG=\pmatrix U&*\\*&*\endpmatrix$ with $U$ an $a\times a$ matrix,
we have
$$\e\left\{W\,^tG\pmatrix I_a\\&0\endpmatrix G/p\right\}
=\e\left\{GW\,^tG\pmatrix I_a\\&0\endpmatrix/p\right\}=\e\{U/p\}.$$
Given an $a\times a$ symmetric matrix $U$ over $\F$, the number of $G\in GL_d(\F)$
so that $GW\,^tG=\pmatrix U&*\\*&*\endpmatrix$ is
$$r^*(W,U)\,p^{a(d-a)}\prod_{i=0}^{d-a-1}(p^{d-a}-p^i)$$
as $\prod_{i=a}^{d-1}(p^d-p^i)$ is the number of ways to extend a $d\times a$ matrix
with rank $a$ over $\F$ to an element of $GL_d(\F)$.  Hence
$$\align
\alpha(W)-1
&=\sum_{a=1}^d\sum_{{U\in\F^{a,a}_{\sym}}\atop{\dim U=a}}
r^*(W,U)\left(\frac{\e\{UI_a/p\}}{o(I_a)}+\frac{\e\{UJ_a/p\}}{o(J_a)}\right)\\
&=\sum_{a=1}^d\sum_{{\cls U}\atop{\dim U=a}}r^*(W,U)
\sum_{G\in O(U)\backslash GL_a}
\left(\frac{\e\{GU\,^tGI_a/p\}}{o(I_a)}
+\frac{\e\{GU\,^tGJ_a/p\}}{o(J_a)}\right)\\
&=\sum_{a=1}^d\sum_{{\cls U}\atop{\dim U=a}}R^*(W,U)
\sum_{G\in GL_a}
\left(\frac{\e\{U\,^tGI_aG/p\}}{o(I_a)}
+\frac{\e\{U\,^tGJ_aG/p\}}{o(J_a)}\right)\\
&=\sum_{a=1}^d\sum_{{\cls U}\atop{\dim U=a}}R^*(W,U)\,\alpha'(U).
\endalign$$
Hence $\alpha(W)=\sum_{a=0}^d R^*(W,U)\alpha'(U)$, as $R^*(W,\{0\})\alpha'(\{0\})=1$.

When $p=2$, the argument is virtually identical.  The difference is that the
 symmetric, invertible $d\times d$ matrices over
$\Z/2\Z$ lie in one $GL_d(\Z/2\Z)$-orbit if $d$ is odd, and in two $GL_d(\Z/2\Z)$-orbits
when $d$ is even; for $d=2c\ge 2$, the $GL_d(\Z/2\Z)$-orbits are represented
by $I_{2c}$ and $J_{2c}=\pmatrix 0&1\\1&0\endpmatrix \perp\cdots\perp
\pmatrix 0&1\\1&0\endpmatrix.$
$\square$
\enddemo

Notice that the formulas in Theorem 3.4 and Corollary 7.3 are identical when the
lattices $\Omega$ and $\Delta$ in Corollary 7.3 are orthogonal.  This leads us to
our final result.

\proclaim{Theorem 7.4}  
Let $f$ be a degree $n$ Siegel modular form of weight $k+1/2$,
level $N$ with $4|N$, and character $\chi$.  Let $\theta^{(n,1)}$ be the 
weight $1/2$, level 4 Jacobi modular form
defined by
$$\theta^{(n,1)}(\tau,Z)=\sum_{R\in\Z^{1,n}} \e\{2\,^tRR\tau+4\,^tRZ\}$$
where $\tau\in\H_{(n)}$ and $Z\in\C^{1,n}$.  Then 
$f(\tau)\mapsto f(\tau)\cdot\theta^{(n,1)}(\tau,Z)$ is an isomorphism
from $\M_{k+1/2}(\Gamma_0^{(n)}(N),\chi)$ onto 
$\J_{k+1,1}^{\text{even}}(\Gamma_0^{(n,1)}(N),\chi')$ 
where $$\chi'(d)=\chi(d)\left({(-1)^{k+1}\over |d|}\right)(\sgn\, d)^{k+1}.$$
Further, for $p$ prime and $1\le j\le n$, we have
$$\big(f|\widetilde T_j(p^2)\big)\cdot\theta^{(n,1)} \,=\,
\big(f\cdot\theta^{(n,1)}\big)|\widetilde T^J_j(p^2)$$
if $p\nmid N$,
$$p^{j/2}\cdot\big(f|T_j(p^2)\big)\cdot\theta^{(n,1)} \,=\,
\big(f\cdot\theta^{(n,1)}\big)|T^J_j(p^2)$$
if $p|N$ and $p\not=2$,
and
$$2^{j/2}\cdot\big(f| T_j(4)\big)\cdot\theta^{(n,1)} \,=\,
\big(f\cdot\theta^{(n,1)}\big)| T^J_j(4)|\psi.$$
\endproclaim

\noindent{\bf Remark.}  When $p\nmid N$, the definitions of
$\widetilde T_j(p^2)$ and $\widetilde T^J_j(p^2)$ give us
$$p^{j/2}\cdot\left(f| T_j(p^2)\right)\theta^{(n,1)}=\left(f\theta^{(n,1)}\right)|T_j(p^2).$$

\smallskip

\demo{Proof}
We first explain why the character of $f(\tau)\,\theta^{(n,1)}(\tau,Z)$
is $\chi'$.
Recall that $\Gamma_0^{(n,1)}(N)$ is generated by matrices of the form
$$\underline\gamma=\pmatrix A&&B\\&1&&0\\C&&D\\&0&&1\endpmatrix,
\ \underline\gamma'=\pmatrix
I&&&^tV\\U&1&V&w\\&&I&-^tU\\&&&1\endpmatrix$$
where $\pmatrix A&B\\C&D\endpmatrix\in\Gamma_0^{(n)}(N)$,
$U,V\in\Z^{1,n}$, $w\in\Z$. 
By Theorem 1.2, for $\underline\tau=\pmatrix \tau&^tZ\\Z&\tau'\endpmatrix
\in\h_{(n+1)}$ with $\tau\in\h_{(n)}$, we have
$$\theta^{(n+1)}(\underline\gamma'\underline\tau)=\theta^{(n+1)}(\underline\tau)
\text{ and }
\theta^{(n+1)}(\underline\gamma\underline\tau)=\frac{\theta^{(n)}(\gamma\tau)}{\theta^{(n)}(\tau)}
\theta^{(n+1)}(\underline\tau).$$
Since $\theta^{(n,1)}(\tau,Z)$ is the 1st Fourier-Jacobi coefficient
of $\theta^{(n+1)}(\underline\tau)$, with
$\widehat\theta^{(n,1)}(\underline\tau)=\theta^{(n,1)}(\tau,Z)\e\{2\tau'\}$
($\tau'$ a formal variable), we have
$$f(\tau)\widehat\theta^{(n,1)}(\tau,Z)|\underline\gamma'
=f(\tau)\widehat\theta^{(n,1)}(\tau,Z),$$
and
$$\align
f(\tau)\widehat\theta^{(n,1)}(\underline\tau)|\underline\gamma
&= (\det(C\tau+D))^{-(k+1)} f(\gamma\tau) 
\widehat\theta^{(n,1)}(\underline\gamma\underline\tau) \\
&= (\det(C\tau+D))^{-(k+1)} 
\left(\frac{\theta^{(n)}(\gamma\tau)}{\theta^{(n)}(\tau)}\right)^{2k+2}
\chi(\det D) f(\tau) \widehat\theta^{(n,1)}(\underline \tau).
\endalign$$
Now, $(\theta^{(n)}(\tau))^{2k+2}=\theta^{(n)}(L;\tau)$ where $L$ is
a lattice with rank $2k+2$ and quadratic form given by the matrix $2I_{2k+2}$.
Hence $\theta^{(n)}(L;\tau)$ is a modular form of weight $k+1$, level 4,
and character $\varphi$ defined by
$$\varphi(d)=\left(\frac{(-1)^{k+1}}{|d|}\right) (\sgn d)^{k+1}.$$
Thus $\theta^{(n)}(L;\gamma\tau)=(\det(C\tau+D))^{k+1} \varphi(\det D)
\theta^{(n)}(L;\tau),$ and hence
$$f(\tau)\widehat\theta^{(n,1)}(\underline\tau)|\underline\gamma
=\chi'(\det D) f(\tau) \widehat\theta^{(n,1)}(\underline\tau).$$
So $f\mapsto f\theta^{(n,1)}$ maps $\M_{k+1/2}(\Gamma_0^{(n)}(N),\chi)$
into $J^{\even}_{k+1,1}(\Gamma^{(n,1)}(N),\chi')$.

Now suppose $F\in J^{\even}_{k+1,1}(\Gamma_0^{(n,1)}(N),\chi').$
Then
$$\align
F(\tau,Z)
&=\sum_{T',R} c(T',2R) \e\{T'\tau+4\,^tRZ\}\\
&=\sum_{T,R} c(T,0) \e\{(T+2\,^tRR)\tau+4\,^tRZ\}\\
&=f(\tau) \theta^{(n,1)}(\tau,Z)
\endalign$$
where $f(\tau)=\sum_{T} c(T,0) \e\{T\tau\}.$
Also, as a function of $Z$, $f(\tau)$ is the 0th Fourier coefficient
of $F$; consequently $f$ is analytic.  Then, reversing the derivation
in the previous paragraph, we see $f\in \M_{k+1/2}(\Gamma_0^{(n)}(N),\chi).$

It is clear that these maps between $\M_{k+1/2}(\Gamma_0^{(n)}(N),\chi)$
and $J^{\even}_{k+1,1}(\Gamma^{(n,1)}(N),\chi')$ are inverses of each other,
and hence we have a bijection between these spaces.

Now take $p$ prime; suppose first that $p\nmid N$ (so in particular, $p\not=2$).
For $f\in \M_{k+1/2}(\Gamma_0^{(n)}(N),\chi)$, set $F=f\theta^{(n,1)}$.  Then
with $\Delta=\Z w\simeq(2)$, we have
$$F(\tau,Z)\,\e\{2\tau'\}
=\sum_{\Lambda}
c_{\Delta}(\Lambda\perp\Delta)\,\e^*_{\Delta}\{(\Lambda\perp\Delta)\underline\tau\}$$
where $c_{\Delta}(\Lambda\perp\Delta)=c_f(\Lambda).$
We first show $F|\widetilde T^J_j(p^2)$ is even.
As a Fourier series, $F(\tau,Z)|\widetilde T^J_j(p^2)\,\e\{2\tau'\}$ is a sum with
exponentials
$$\e\{\underline T[X_j^{-1}D\underline G^{-1}]\underline\tau\}$$
where 
$X_j=\diag\{pI_j,I_{n+1-j}\}$, $D=\diag\{I_{r_0},pI_{r_1},p^2I_{r_2},I_{n+1-j}\}$
with $r_0+r_1+r_2=j$, $\underline G\in GL_{n+1,1}^J(\Z)$,  
$\underline T=\pmatrix T+2\,^tRR&2\,^tR\\2R&2\endpmatrix$ and
$\underline T[X_j^{-1}D\underline G^{-1}]$
even integral.
Consequently 
$$\underline T[X_j^{-1}D\underline G^{-1}]
=\pmatrix T'&2\,^tR'\\2R'&2\endpmatrix
\text{ where }
R'\in {1\over p}\Z^{1,n}\cap{1\over 2}\Z^{1,n}=\Z^{1,n}.$$
Thus with $\Lambda'\oplus\Delta\simeq\pmatrix T'&2\,^tR'\\2R'&2\endpmatrix,$
we have
$$\Lambda'\oplus\Delta = (\Lambda'\oplus\Delta)\pmatrix I&0\\-R'&1\endpmatrix
\simeq\pmatrix T'-2\,^tR'R'&0\\0&2\endpmatrix.$$
Hence 
$$\big(F(\tau,Z)|\widetilde T^J_j(p^2)\big)\,\e\{2\tau'\}
=\sum_{\Lambda\perp\Delta} \widetilde c_{\Delta}(\Lambda\perp\Delta)
\,\e^*_{\Delta}\{(\Lambda\perp\Delta)\underline\tau\};$$
by Theorem 7.2,
$$\widetilde c_{\Delta}(\Lambda\perp\Delta)
=\sum_{\Omega'}\widetilde A^J_{j,\Delta}(\Lambda,\Omega') c_{\Delta}(\Omega'\oplus\Delta)$$
where $\Omega'$ varies so that $
p\Lambda\perp\Delta\subseteq\Omega'\oplus\Delta\subseteq{1\over p}(\Lambda\perp\Delta)$
and $p\Lambda'\subseteq\Omega'\subseteq{1\over p}\Lambda'$ for some $\Lambda'$
with $\Lambda'\oplus\Delta=\Lambda\perp\Delta$.
(Recall that $c_{\Delta}(\Omega'\oplus\Delta)=0$ unless $\Omega'\oplus\Delta$ is
even integral.)  
Since $\Delta\simeq(2)$ and $p\not=2$, any even integral sublattice of 
${1\over p}(\Lambda\perp\Delta)$ is actually a sublattice of
${1\over p}\Lambda\perp\Delta$.  So suppose
$$p\Lambda\perp\Delta\subseteq\Omega'\oplus\Delta\subseteq
{1\over p}\Lambda\perp\Delta.$$
Taking $\Omega=(\Omega'\oplus\Delta)\cap{1\over p}\Lambda$, we have
$$\Omega\perp\Delta=\Omega\oplus\Delta=
(\Omega'\oplus\Delta)\cap\left({1\over p}\Lambda\perp\Delta\right)=\Omega'\oplus\Delta.$$
Consequently
$$\align
\widetilde c_{\Delta}(\Lambda\perp\Delta)
&=\sum_{p\Lambda\subseteq\Omega\subseteq{1\over p}\Lambda}
\widetilde A^J_{j,\Delta}(\Lambda,\Omega) c_{\Delta}(\Omega\perp\Delta)\\
&=\sum_{p\Lambda\subseteq\Omega\subseteq{1\over p}\Lambda}
\widetilde A_j(\Lambda,\Omega) c_f(\Omega),
\endalign$$
which is the $\Lambda$th coefficient of $f|\widetilde T_j(p^2)$, and hence
the $(\Lambda\perp\Delta)$th coefficient of
$\left(f(\tau)|\widetilde T_j(p^2)\right)\cdot\theta^{(n,1)}(\tau,Z).$

Now suppose $p|N$.  When $p\not=2$, the argument follows much as when
$p\nmid N$, except the situation is simpler since $D$ is always $I$.  When $p=2$,
this argument breaks down since $\frac{1}{p}\Z^{1,n}\cap\frac{1}{2}\Z^{1,n}
=\frac{1}{2}\Z^{1,n}$, not $\Z^{1,n}$.  For this reason we need to follow
$T^J_j(4)$ by $\psi$ to get the desired equality. $\square$
\enddemo

\enddocument